\documentclass[letterpaper, 11pt]{article}
\usepackage{fullpage}
\usepackage{amssymb,amsmath}
\usepackage{epsf}
\usepackage{times,bm,mathrsfs}  % Fonts
\usepackage{amsmath}
\usepackage{amssymb}
\usepackage{amsfonts}
\usepackage{mathrsfs}
\usepackage{bbm}
\usepackage{color}
\usepackage{graphicx}
\usepackage{amscd}
\usepackage{epsfig}

\usepackage{geometry}
\definecolor{Red}{rgb}{1,0,0}
\definecolor{Blue}{rgb}{0,0,1}
\definecolor{Olive}{rgb}{0.41,0.55,0.13}
\definecolor{Green}{rgb}{0,1,0}
\definecolor{MGreen}{rgb}{0,0.8,0}
\definecolor{DGreen}{rgb}{0,0.55,0}
\definecolor{Yellow}{rgb}{1,1,0}
\definecolor{Cyan}{rgb}{0,1,1}
\definecolor{Magenta}{rgb}{1,0,1}
\definecolor{Orange}{rgb}{1,.5,0}
\definecolor{Violet}{rgb}{.5,0,.5}
\definecolor{Purple}{rgb}{.75,0,.25}
\definecolor{Brown}{rgb}{.75,.5,.25}
\definecolor{Grey}{rgb}{.5,.5,.5}
\definecolor{Black}{rgb}{0,0,0}

\newcommand{\acal}{\mathcal{A}}
\newcommand{\bcal}{\mathcal{B}}
\newcommand{\ccal}{\mathcal{C}}

\newcommand{\ecal}{\mathcal{E}}

\newcommand{\lcal}{\mathcal{L}}

\newcommand{\scal}{\mathcal{S}}

\newcommand{\nintgr}{\mathbb{N}}

\newcommand{\eps}{\varepsilon}

\newcommand{\ind}{\mathbbm{1}}

\newcommand{\bdm}{\begin{displaymath}}
\newcommand{\edm}{\end{displaymath}}
\newcommand{\bea}{\begin{eqnarray*}}
\newcommand{\eea}{\end{eqnarray*}}
\newcommand{\bean}{\begin{eqnarray}}
\newcommand{\eean}{\end{eqnarray}}

\newcommand{\prob}{\mathbb{P}}
\newcommand{\expec}{\mathbb{E}}

\newcommand{\poly}{\mathrm{poly}}

\newtheorem{theorem}{Theorem}[section]
\newtheorem{lemma}[theorem]{Lemma}

\newtheorem{remark}[theorem]{Remark}

\newtheorem{definition}[theorem]{Definition}

\newenvironment{app-proof}[1]{\noindent{\textbf{Proof of #1:}}}{$\blacksquare$\vskip\belowdisplayskip}

\newcommand{\tree}[2]{T^{(#1)}_{#2}}
\newcommand{\prbs}{p_{\mathrm{s}}}
\newcommand{\ths}{\theta_{\mathrm{s}}}
\newcommand{\prbid}{p_{\mathrm{id}}}
\newcommand{\prbi}{p_{\mathrm{i}}}
\newcommand{\prbd}{p_{\mathrm{d}}}
\newcommand{\arm}[2]{R^{(#1)}_{#2}}
\newcommand{\bern}[2]{B_{#1,#2}}

\newcommand{\Maj}{\mathrm{Maj}}

\newcommand{\x}{x}
\newcommand{\X}{X}
\newcommand{\rx}{\hat \x}
\newcommand{\RX}{\widehat \X}
\newcommand{\y}{y}
\newcommand{\Y}{Y}
\newcommand{\ry}{\hat \y}
\newcommand{\RY}{\widehat \Y}

\newcommand{\kref}{k}
\newcommand{\rkref}{\hat \kref}
\newcommand{\nisle}{N}
\newcommand{\lisle}{\ell}
\newcommand{\rs}{\hat s}
\newcommand{\thres}{\gamma}
\newcommand{\nref}{n}

\newcommand{\mink}{\underline{\kref}}
\newcommand{\maxk}{\bar \kref}
\newcommand{\maxnisle}{\overline{\nisle}}

\newcommand{\path}{\mathrm{Path}}
\newcommand{\RXi}{\widehat{\Xi}}
\newcommand{\rxi}{\hat{\xi}}
\newcommand{\Href}{H}

\newcommand{\bin}{\mathrm{Bin}}

\newcommand{\spin}[1]{\langle #1 \rangle}

\newcommand{\RA}{\widehat{A}}
\newcommand{\RD}{\widehat{D}}
\newcommand{\RI}{\widehat{I}}
\renewcommand{\AA}{\mathscr{A}}
\newcommand{\DD}{\mathscr{D}}
\newcommand{\II}{\mathscr{I}}
\newcommand{\RAA}{\hat{\AA}}
\newcommand{\RDD}{\hat{\DD}}
\newcommand{\RII}{\hat{\II}}

\DeclareMathOperator{\cor}{Corr}

\author{
Alexandr Andoni\footnote{CSAIL, MIT}, Constantinos Daskalakis\footnote{CSAIL, MIT. {\tt costis@mit.edu}. Part of this work was done while the author was a postdoctoral researcher at Microsoft Research.},\\ Avinatan Hassidim\footnote{MIT}, Sebastien Roch\footnote{Department of Mathematics, UCLA. Part of this work was done while the author was a postdoctoral researcher at Microsoft Research.}
}

\title{\vspace{0cm}
Global Alignment of Molecular Sequences\\ via Ancestral State Reconstruction}
\begin{document}

\maketitle
\addtocounter{page}{-1}
\thispagestyle{empty}

\begin{abstract}
Molecular phylogenetic techniques do not generally account for such
common evolutionary events as site insertions and deletions (known as
indels). Instead tree building algorithms and ancestral state
inference procedures typically rely on substitution-only models of
sequence evolution. In practice these methods are extended beyond this
simplified setting with the use of heuristics that produce global
alignments of the input sequences---an important problem which has no
rigorous model-based solution. In this paper we open a new direction
on this topic by considering a version of the multiple sequence
alignment in the context of stochastic indel models. More precisely,
we introduce the following {\em trace reconstruction problem on a
tree} (TRPT): a binary sequence is broadcast through a tree channel
where
we allow substitutions, deletions, and insertions; we seek
to reconstruct the original sequence from the sequences
received at the leaves of the tree. We give a recursive
procedure for this problem with strong reconstruction
guarantees at low mutation rates, providing also an alignment of the
sequences at the leaves of the tree. The TRPT problem without indels
has been studied in previous work (Mossel 2004, Daskalakis et al.
2006) as a bootstrapping step towards obtaining
information-theoretically optimal phylogenetic reconstruction methods.
The present work sets up a framework for extending these works to
evolutionary models with indels.

In the TRPT problem we begin with a random sequence $x_1, \ldots, x_k$ at the root of a $d$-ary tree. If vertex $v$ has the sequence $y_1, \ldots y_{k_{v}}$, then each one of its $d$ children will have a sequence which is generated from $y_1, \ldots y_{k_{v}}$ by flipping three biased coins for each bit. The first coin has probability $p_s$ for Heads, and determines whether this bit will be substituted or not. The second coin has probability $p_d$, and determines whether this bit will be deleted, and the third coin has probability $p_i$ and determines whether a new random bit will be inserted. The input to the procedure is the sequences of the $n$ leaves of the tree, as well as the tree structure (but not the sequences of the inner vertices) and the goal is to reconstruct an approximation to the sequence of the root (the DNA of the ancestral father). For every $\epsilon > 0$,  we present a deterministic algorithm which outputs an approximation of $x_1, \ldots, x_k$ if $\prbi 
 + \prbd < O(1/k^{2/3} \log n)$ and $(1 - 2\prbs)^2 > O(d^{-1} \log d)$. 

%We argue that our requirement on $p_s$ is close to optimal, and that our bounds on $p_i$, $p_d$ have the right dependence on $k$ and $n$ to within polynomial factors.

To our knowledge, this is the first rigorous trace reconstruction result on a tree in the presence of indels.
\end{abstract}

\clearpage

\section{Introduction}

\noindent\textbf{Trace reconstruction on a star.}
In the ``trace reconstruction
problem'' (TRP)~\cite{Levenshtein:01a, Levenshtein:01b,
BaKaKhMG:04, KannanMcGregor:05, HoMiPaWi:08, ViswanathanSwaminathan:08},
a random binary string
$\X$ of length $\kref$
generates an i.i.d.~collection of traces $\Y_1,\ldots,\Y_\nref$ that
are identical to $\X$ except for random \emph{mutations} which
consist in \emph{indels}, i.e., the deletion of an old site or the insertion of
a new site between existing sites, and \emph{substitutions},
i.e., the flipping of the state at an existing
site\footnote{One can also consider the
case where $\X$ is arbitrary rather than random. We will not discuss
this problem here.}.
(In keeping with biological terminology, we refer to the
components or positions of a string as \emph{sites}.)
The goal is to reconstruct efficiently the original
string with high probability from as few random traces as possible.

An important motivation for this problem
is the reconstruction of ancestral DNA sequences in computation biology~\cite{BaKaKhMG:04, KannanMcGregor:05}.
One can think of $\X$ as a gene in an (extinct) ancestor species $0$.
Through speciation, the ancestor $0$ gives rise to a large
number of descendants $1,\ldots,\nref$ and gene $\X$
evolves independently through the action of mutations into
sequences $\Y_1,\ldots,\Y_\nref$ respectively.
Inferring the sequence $\X$ of an ancient gene from extant descendant copies
$\Y_1,\ldots,\Y_\nref$ is a standard problem in evolutionary
biology~\cite{Thornton:04}.
The inference of $\X$ typically
requires the solution of an auxiliary problem,
the \emph{multiple sequence alignment problem} (which is
an important problem in its own right in computational
biology): site $t_i$ of sequence $\Y_i$ and
site $t_j$ of sequence $\Y_j$ are said to be \emph{homologous} (in this
simplified TRP setting) if
they descend from a common site $t$ of $\X$ \emph{only through
substitutions}; in the multiple sequence alignment
problem, we seek roughly to uncover the homology relation
between $\Y_1,\ldots,\Y_\nref$.
Once homologous sites have been identified, the original sequence
$\X$ can be estimated, for instance, by site-wise majority.

The TRP as defined above is an \emph{idealized}
version of the ancestral sequence reconstruction problem in
one important aspect. It ignores the actual phylogenetic
relationship between species $1,\ldots,\nref$.
A \emph{phylogeny} is a (typically, binary) tree relating a group of species.
The leaves of the tree correspond to extant species.
Internal nodes can be thought of as extinct ancestors.
In particular the root of the tree represents the most
recent common ancestor of all species in the tree.
Following paths from the root to the leaves, each bifurcation indicates
a speciation event whereby two new species are created
from a parent. An excellent introduction to phylogenetics is \cite{SempleSteel:03}.

A standard assumption in computational phylogenetics
is that genetic information evolves from
the root to the leaves according to a Markov model on the tree.
%It is further assumed that this process is repeated independently a number of times
%denoted $k$.
%Thus each node of the tree is associated with a sequence of length $k$.
%The vector of the $i$'th letter of all sequences at the leaves is called the
%$i$'th {\em character}. One of the major tasks in
%molecular biology, the {\em reconstruction of phylogenetic trees},
%is to infer the topology of the tree from the characters
%at the leaves.
Hence, the stochastic model used in
trace reconstruction can be seen as a special case
where the phylogeny is \emph{star-shaped}.
(The substitution model used in trace reconstruction
is known in biology as the Cavender-Farris-Neyman
(CFN)~\cite{Cavender:78,Farris:73,Neyman:71} model.)
It may seem that a star is a good first approximation for the
evolution of DNA sequences. However extensive work
on the so-called ``reconstruction problem'' in theoretical
computer science and statistical physics has
highlighted the importance of taking into account the full tree model
in analyzing the reconstruction of ancestral sequences.
%In this paper we will be mostly interested in two evolutionary models,
%the so-called  and Jukes-Cantor
%(JC)~\cite{JukesCantor:69} models.
%In the CFN model the states at the nodes of the tree
%are $0$ and $1$ and their a
%priori probability at the root is uniform. To
%each edge $e$ corresponds a {\em mutation probability} $p(e)$ which is the
%probability that the state changes along the edge $e$.
%Note that this model is identical to the free Gibbs measure of the Ising model on the tree.
%See~\cite{Lyons:89}.
%In the JC model the states are $A$, $C$, $G$ and $T$ with a priori
%probability $1/4$ each. To each edge $e$ corresponds a mutation
%probability $p(e)$ and it is assumed that every state transitions with
%probability $p(e)$ to each of the other states.
%This model is equivalent to the ferromagnetic Potts model on the tree.

\noindent\textbf{The ``reconstruction problem.''}
In the ``reconstruction problem'' (RP),
we have a single site which evolves through substitutions
only from the root to the leaves of a tree.
In the most basic setup which we will consider here,
the tree is $d$-ary and each edge is an independent symmetric indel-free channel where
the probability of a substitution is a constant
$\prbs > 0$. The goal is to reconstruct the state at the
root given the vector of states at the leaves.
More generally, one can consider a sequence of length $\kref$ at the root
where each site evolves independently according to the Markov process above.
Denote by $\nref$ the number of leaves in the tree.
The RP has attracted much attention in the theoretical
computer science literature due to its deep connections to
computational phylogenetics~\cite{Mossel:03,
Mossel:04a,DaMoRo:06,Roch:08}
and statistical physics~\cite{Mossel:98,EvKePeSc:00,Mossel:01,MosselPeres:03,
MaSiWe:04,JansonMossel:04,BeKeMoPe:05,BoChMoRo:06,
GerschenfeldMontanari:07,BhVeViWe:arxiv,Sly:09a,Sly:09b}.
See e.g.~\cite{Roch:thesis,Sly:thesis} for background.

Unlike the star case, the RP on a tree exhibits
an interesting thresholding effect:
%\textbf{(Information loss)}
on the one hand,
information is lost at an
exponential rate along each path from the root;
%\textbf{(Exponential growth)}
on the other hand, the number of
paths grows exponentially with the number of levels.
When the substitution probability is low, the latter ``wins''
and vice versa.
This ``phase transition'' has been thoroughly analyzed in the theoretical
computer science and mathematical physics literature---although
much remains to be understood.
More formally, we say that the RP is {\em solvable}
when the correlation between the root and the leaves persists no
matter how large the tree is. Note that unlike the TRP we do not require
high-probability reconstruction in this case as it is
not information-theoretically achievable for $d$ constant---simply
consider the information lost on the
first level below the root. Moreover the ``number of traces''
is irrelevant here as it is governed by the depth of the tree and
the solvability notion implies nontrivial correlation
for any depth. When the RP is unsolvable, the
correlation decays to $0$ for large trees. The results
of~\cite{BlRuZa:95,EvKePeSc:00,Ioffe:96a,BeKeMoPe:05,MaSiWe:04,BoChMoRo:06}
show that for the CFN model, if $\prbs < p^{\ast}$,
then the RP is
solvable, where
$d (1 - 2p^{\ast})^2 = 1$.
This is the so-called \emph{Kesten-Stigum} bound~\cite{KestenStigum:66}.
If, on the other hand,
$\prbs > p^{\ast}$, then the
RP is {\em unsolvable}. Moreover in this case, the
correlation between the root state and any function of the character
states at the leaves decays as $\nref^{-\Omega(1)}$.
The positive result above is obtained by taking a majority
vote over the leaf states.

Like the TRP, the RP is only an \emph{idealized} version of the
ancestral sequence reconstruction problem: it ignores
the presence of indels. In other words, the RP assumes that
the multiple sequence alignment problem has been solved perfectly.
This is in fact a long-standing assumption in evolutionary biology
where one typically preprocesses sequence data by running it
through a multiple sequence alignment heuristic and then one
only has to model the substitution process.
This simplification has come under attack in the
biology literature, where it has been argued that
alignment procedures often create systematic
biases that affect analysis~\cite{LoytynojaGoldman:08, WoSuHu:08}.
Much empirical work has been devoted
to the proper joint estimation of alignments
and phylogenies~\cite{ThKiFe:91,ThKiFe:92,Metzler:03,
MiLuHo:04,SuchardRedelings:06,RivasEddy:08,
LoytynojaGoldman:08, LRNLW:09}.

\noindent\textbf{Our results.} 
We make progress in this recent new direction by analyzing the RP in
the presence of indels---which we also refer to as the TRP on a tree
(TRPT).  We consider a $d$-ary tree where each edge is an independent
channel
with substitution probability $\prbs$, deletion probability $\prbd$,
and insertion probability $\prbi$ (see Section~\ref{sec:model} for a
precise statement of the
model). The root sequence has length $\kref$ and is assumed to be
uniform in $\{0,1\}^\kref$. As in the standard RP, we drop the requirement
of high-probability reconstruction and seek instead a reconstructed sequence
that has correlation with the true root sequence uniformly bounded in the depth.

We give an efficient recursive procedure which solves the TRPT
for $\prbs > 0$ a small enough constant (strictly below, albeit close,
to the Kesten-Stigum bound)
and $\prbd,\prbi = O(\kref^{-2/3}\log^{-1} \nref)$.
As a by-product of our analysis we also obtain
a partial global alignment of the sequences at the leaves. Our method
provides a framework for separating the indel process from the
substitution process by identifying well-preserved subsequences which
then serve as markers for alignment and reconstruction (see
Section~\ref{sec:results} for a high-level description of our
techniques). As far as we are aware, our results are the first
rigorous results for this problem.

Results on the RP have been used in previous work to advance the state
of the art in rigorous phylogenetic tree reconstruction
methods~\cite{Mossel:04a,DaMoRo:06,MiHiRa:08,Roch:08}. A central
component in these methods is to solve the RP on a partially
reconstructed phylogeny to obtain sequence information that is
``close'' to the evolutionary past; then this sequence information is
used to obtain further structural information about the phylogeny. The
whole phylogeny is built by alternating these steps. Our method sets
up a framework for extending these techniques beyond substitution-only
models. Partial results of this type will be given in the full version
of the paper.

\noindent\textbf{Related work.}
Much work has been devoted to the
trace reconstruction problem on a star~\cite{Levenshtein:01a, Levenshtein:01b,
BaKaKhMG:04, KannanMcGregor:05, HoMiPaWi:08, ViswanathanSwaminathan:08}.
In particular, in~\cite{HoMiPaWi:08},
it was shown that, when there are only deletions, it is possible to tolerate a small constant deletion rate using  $\poly(\kref)$ traces. For a different range of parameters,
Viswanathan and Swaminathan~\cite{ViswanathanSwaminathan:08}
showed that, under constant substitution probability and
$O(1/\log \kref)$ indel probability, $O(\log \kref)$ traces suffice.
Both results assume that the root sequence $\X$ is uniformly random.

The multiple sequence alignment problem as
a combinatorial optimization problem (finding the best
alignment under some pairwise scoring function)
is known to be NP-hard~\cite{WangJiang:94,Elias:06}.
Most heuristics used in practice, such as
CLUSTAL~\cite{HigginsSharp:88},
T-Coffee~\cite{NoHiHe:00},
MAFFT~\cite{KaMiKuMi:02},
and
MUSCLE~\cite{Edgar:04},
use the idea of a guide tree, that is,
they first construct a very rough phylogenetic
tree from the data (using edit distance
as a measure of evolutionary distance),
and then recursively construct
local alignments produced by ``aligning alignments.''
Our work can be thought of as the first attempt to analyze
rigorously this type of procedure.

%In a recent related work, two of us have also
%given the first efficient phylogenetic reconstruction
%algorithm in the presence of indels~\cite{DaskalakisRoch:09}.

Finally, our work is tangentially related to the study of edit distance.
Edit distance and pattern matching in random environments
has been studied, e.g., by \cite{navarro2001gta,navarro:ima,andoni2008sce}.
%It is known that in these cases one can obtain much stronger results than for
%worst-case instances.
%Moreover, it is known that practical heuristics
%perform better than the guarantees given by the smoothed analysis,
%or the random case analysis, and it is an interesting question what
%is the exact property which enables heuristics to perform better
%than the proven theoretical results.

%\srnote{Alex, can you add some references to the edit distance
%literature? Comment on how the problems ``differ.''}

\subsection{Definitions} \label{sec:model}

We now define our basic model of sequence evolution.
\begin{definition}[Model of sequence evolution] \label{def: model of evolution}
Let $\tree{d}{H}$ be the $d$-ary tree with $H$ levels and $n=d^H$ leaves.
For simplicity, we assume throughout that $d$ is odd.
We consider the following model of evolution on $\tree{d}{\Href}$.
The sequence at the root of $\tree{d}{\Href}$ has length $\kref$ and is drawn uniformly at random
over $\{0,1\}^{\kref}$.
Along each edge of the tree, each site (or position) undergoes the following
mutations independently of the other sites:
\begin{itemize}
\item
\textbf{Substitution.} The site state is flipped with
probability $\prbs > 0$.
\item
\textbf{Deletion.}
The site is deleted with probability $\prbd > 0$.
\item
\textbf{Insertion.}
A new site is created to the right of the current site with probability $\prbi > 0$.
The state of this
new site is uniform $\{0,1\}$.
\end{itemize}
These operations occur independently of each other.
The last two are called \emph{indels}.
We let $\prbid = \prbi + \prbd$
and $\ths = 1-2 \prbs$. The parameters $\prbs, \prbd, \prbi$
may depend on $\kref$ and $\nref$,
where $\nref$ is the number of leaves.
\end{definition}
\begin{remark}
For convenience, our model of insertion is intentionally simplistic.
In the biology literature, related continuous-time Markov models
are instead used for this kind of process~\cite{ThKiFe:91,ThKiFe:92,Metzler:03,
MiLuHo:04,RivasEddy:08,DaskalakisRoch:09}.
It should be possible to extend our results to such
generalizations by proper modifications to the algorithm.
\end{remark}

\subsection{Results}\label{sec:results}

\noindent\textbf{Statement of results.}
Our main result is the following.
Denote by $\X = \x_1,\ldots,\x_\kref$ a binary
uniform sequence of length $\kref$.
Run the evolutionary process on $\tree{d}{H}$
with root sequence $\X$
and let $\Y_1,\ldots,\Y_\nref$ be the sequences
obtained at the leaves, where
$\Y_i = \y^i_1,\ldots,\y^i_{\kref_i}$.
\begin{theorem}[Main result] \label{thm:main}
For all $\chi > 0$ and
$\beta = O(d^{-1} \log d)$, there is
$\Phi, \Phi', \Phi'' > 0$ such that the following holds
for $d$ large enough.
There is a polynomial-time
algorithm $\mathbb{A}$ with access to $\Y_1,\ldots,\Y_\nref$
such that for all
\begin{equation*}
(1-2\prbs)^2 >  \frac{\Phi \log d}{d},
%(1-2\prbs)^2 > \max\left\{\frac{\Phi \log d}{d}, 16 \beta\right\},
%\end{equation*}
%\begin{equation*}
\qquad\prbi + \prbd < \frac{\Phi'}{\kref^{2/3} \log \nref},
%\end{equation*}
%and
%\begin{equation*}
\qquad \Phi''\log^3 \nref < \kref < \poly(n),
\end{equation*}
the algorithm $\mathbb{A}$ outputs a binary sequence $\RX$ which satisfies
the following with probability at least $1 -\chi$:
\begin{enumerate}
\item $\RX = \rx_1,\ldots,\rx_\kref$ has length $\kref$.
\item For all $j = 1,\ldots,\kref$,
$\prob[\rx_j = \x_j] > 1 - \beta$.
\end{enumerate}
\end{theorem}
\begin{remark}
Notice that we assume that the (leaf-labelled) tree and
and the sequence length of the root are known. 
The requirement that the sequence length is known is not crucial. 
We adopt it for simplicity in the presentation.
\end{remark}
\begin{remark}
In fact, we prove a stronger result which allows
$\chi = o(1)$ and shows that the ``agreement'' between $\RX$
and $\X$ ``dominates'' an i.i.d.~sequence.
See Lemma~\ref{lem:limit} and Section~\ref{sec:correctness}.
\end{remark}

\smallskip \noindent\textbf{Proof sketch.}
We give a brief proof sketch.
As discussed previously, in the presence of indels
the reconstruction of ancestral sequences requires
the solution of the \emph{multiple sequence alignment}
problem. However, in addition to being computationally
intractable, global alignment through the optimization
of a pairwise scoring function may create biases
and correlations that are hard to quantify.
%For instance, an optimal solution may have a systematic tendency
%to align sites in the same state, potentially amplifying
%the effect of mutations.
Therefore, we require a more probabilistic approach.
From a purely information-theoretic point of view
the pairwise alignment of sequences that are far apart
in the tree is difficult. A natural solution to
this problem is instead to perform \emph{local}
alignments and ancestral reconstructions,
and recurse our way up the tree.

This \emph{recursive} approach raises its own set of issues.
Consider a parent node and its $d$ children. It may be easy
to perform a local alignment of the children's sequences
and derive a good approximation
to the parent sequence (for example, through site-wise majority).
Note however that, to allow
a recursion of this procedure all the way to the root,
we have to provide strong guarantees about the
probabilistic behavior of our local ancestral reconstruction.
As is the case for global alignment, a careless
alignment procedure creates biases and correlations
that are hard to control. For instance, it is tempting to treat misaligned sites
as independent unbiased noise but this idea presents difficulties:
\begin{quote}
Consider a site $j$ of the parent sequence and
suppose that for this site we have succeeded in aligning
all but two of the children, say $1$ and $2$. Let $\x_{j_i}^i$ denote the site in the $i$'th child which was used to estimate the $j$'th site.
By the independence
assumption on the root sequence and the inserted sites,
$\x_{j_1}^1$ and $\x_{j_2}^2$ are uniform and independent of
$(\x_{j_i}^i)_{i=3}^{d}$.
However, $\x_{j_1}^1$ and $\x_{j_2}^2$
may originate from the \emph{same} neighboring site of the parent
sequence and therefore are themselves correlated.
\end{quote}
Quantifying the effect
of this type of correlation appears to be nontrivial.

Instead, we use an \emph{adversarial} approach to local
ancestral reconstruction.
That is, we treat
the misaligned sites as being controlled by an adversary
who seeks to flip the reconstructed value.
This comes at a cost: it produces an
asymmetry in our ancestral reconstruction. Although the RP
is well-studied in the symmetric noise case, much
remains to be understood in the asymmetric case. In particular,
obtaining tight results in terms
of substitution probability here may not be possible as the critical
threshold of the RP may be hard to identify.
We do however provide a tailored analysis of the particular
instance of the RP by recursive majority
obtained through this adversarial
approach and we obtain results that are close to the known
threshold for the symmetric case. Unlike the standard RP,
the reconstruction error is not i.i.d.~but we show instead that it
``dominates'' an i.i.d.~noise. (See Section~\ref{sec:estimation} for a definition.)
This turns out to be enough
for a well-controlled recursion.
{%Handling the errors in the reconstruction is done in two steps.
We first define a local alignment procedure which has a fair success
probability (independent of $\nref$).  However, applying this alignment procedure
multiple times in the tree is bound to fail sometimes.
We therefore prove that the local reconstruction procedure is somewhat robust
in the sense that even if one of the $d$ inputs to the reconstruction
procedure is faulty, it still has a good probability of success.
%We then analyze the global behavior of the procedure.
%Working with two layers of correction, enables us to handle higher probabilities of insertions and deletions.
}

%As for our local alignment procedure, we adopt an \emph{anchor} approach.
%Anchors were also used by~\cite{KannanMcGregor:05, HoMiPaWi:08}---although
%in a quite different way. We divide the parent sequence into
%islands of length $O(\kref^{1/3})$. (The choice
%of island length $k^{1/3}$ comes from a trade-off between
%the length and number of islands in bounding the ``bad'' events below.
%See the proof of Lemma~\ref{lem:badIndelStructure}.)
%At the beginning of each island
%we have an anchor of length $O(\log \nref)$. Taking an indel probability of %$O(1/\kref^{2/3} \log \nref)$, this gives that with $1 - o(1)$ probability %(but not w.h.p.)
%1) There are a polynomial number of deletions in each link;
%2) Each site is deleted in at least one of the paths from the root to the %leaves; 3) The anchors of the father are indel free (because they are short); %And 4) In every parent island, almost all the corresponding children
%islands have no indel at all and, moreover, that at most one child island may %have a single indel.

As for our local alignment procedure, we adopt an \emph{anchor} approach.
Anchors were also used by~\cite{KannanMcGregor:05, HoMiPaWi:08}---although
in a quite different way. We imagine a partition of every node's sequence into
islands of length $O(\kref^{1/3})$. (The precise choice
of the island length comes from a trade-off between
the length and the number of islands in bounding the ``bad'' events below---see the proof of Lemma~\ref{lem:badIndelStructure}.) At the beginning of each island
we have an anchor of length $O(\log \nref)$. Through this partition of the sequences in islands and anchors we aim to guarantee the following. Given a specific father node $v$,
%(without applying a union bound on the entire tree),
with fair probability
1) all the anchors in the children nodes are indel-free; and
2) for all parent islands, almost all of the corresponding children
islands have no indel at all and, moreover, 
at most one child island may have a single indel.
The ``bad'' children islands---those that do not satisfy these
properties---are treated as controlled by an
adversary. We show that Conditions 1) and 2) are sufficient to guarantee that: the anchors of all islands can be aligned with high probability and single indel events
between anchors can be identified. This allows a local alignment
of all islands with at most one ``bad'' child per island
and is enough to perform a successful adversarial recursive majority vote
as described above. The bound on the maximum indel probability sustained by our reconstruction algorithm comes from satisfying Conditions 1) and 2) above. 

\noindent\textbf{Notation.} For a sequence $X = x_1,\ldots,x_k$,
we let $X[i:j] = x_i,\ldots,x_j$. We use the expression
``with high probability (w.h.p.)'' to mean ``with
probability at least $1-1/\poly(\nref)$'' where the
polynomial in $n$ can be made of arbitrarily high degree
(by choosing the appropriate constants large enough).
We denote by $\bin(n,p)$ a random variable with binomial
distribution of parameters $n,p$. For two random variables
$X,Y$ we denote by $X\sim Y$ the equality in distribution.

\noindent\textbf{Organization.}
The rest of the paper is organized as follows.
We describe the algorithm in Section~\ref{sec:algo}.
The proof of our main result is divided into
two sections. In Section~\ref{sec:analysis},
we prove a series of high-probability claims
about the evolutionary process. Then,
conditioning on these claims, we provide
a deterministic analysis of the correctness
of the algorithm in Section~\ref{sec:correctness}.
All proofs are in the Appendix.

\section{Description of the Algorithm}\label{sec:algo}

In this section we describe our algorithm for TRPT.
Our algorithm is recursive, proceeding from the leaves of the tree to the root. We 
describe the recursive step applied to a non-leaf node of
the tree.

{ \noindent\textbf{Recursive Setup---Our Goal.}}
For our discussion in this section, let us consider a non-leaf node $v$ with $d$
children, denoted $u_i$ for $i\in[d]$. For notational convenience, we drop the index $u$ and denote its children by $1,
\ldots, d$. Our goal for the recursive step of the algorithm is to reconstruct the sequence at the node $v$
given the sequences of the children.  Denote the sites of the father
by $\X_0 = \x^0_1, \ldots, \x^0_{\kref_0}$, and the sites of the
$i$'th child by $\X_i = \x^{i}_1, \ldots, \x^{i}_{\kref_i}$.  During the reconstruction process, we do not have access to the children's sequences, but rather to reconstructed sequences denoted by
$\RX_{i}=\rx^{i}_1, \ldots \rx^{i}_{\rkref_i}$.
%In the analysis later we assume
%that the sites of one child are chosen in an adversarial manner, and
%for the rest of the children, the event that $\rx^{i}_t = \x^{i}_{t}$ stochastically
%dominates an independent random event which
%happens with probability $1 - \beta$.

%Let $\kref$ be a constant to be determined later. (We will show that
%$\kref_0 \approx \kref$ with high probability.)
Let us consider the following partition of the sequence of $v$ into subsequences, called {\it islands}. Of course our algorithm doesn't have access to the sequence at $v$ during the recursive step of the algorithm. We define the partition as a means to describe our algorithm: The sites
of $v$ are partitioned into {\it islands} of length $\lisle = \kref^{1/3}$
(except for the last one which is possibly shorter).
%\srnote{What should
%  we do with the last island? (I think the islands should have a
%  deterministic length.)  I'm suggesting just getting rid of
%  them. Over log n levels, that should hurt us too much. May have to
%  change the statement though.  Alternatively, we could reconstruct
%  everything again but starting form the end. Then match the two
%  reconstructions.}
Denote by $\nisle_0 = \lceil\kref_0/\lisle\rceil$ the number of islands in $v$.  Each island
starts with an {\it anchor} of $a$ bits.
%\srnote{I'm starting all the
%  strings at 1.}
That is, the islands are the bitstrings
$\X_0[1:\lisle], \ \X_0[\lisle+1:2\lisle], \ldots$ and the anchors are
the bitstrings $\X_0[1:a], \ \X_0[\lisle+1:\lisle + a],
\ldots$.

Our algorithm tries to identify for each island
$\X_0[(i-1)\lisle+1:i\lisle]$ the substrings of each of the $d$
children that correspond to this island (i.e., contain the sites of
the island), called ``child islands.'' We do so iteratively for
$i=1\ldots N_0$. We use the islands that did not have indels for
sequence reconstruction, using the substitution-only model. Some
islands will have indels however. This leads to two ``modes of
failure'': one invalidates the entire (parent) node, and the other
invalidates only an island of a child. More specifically, a node becomes invalidated
(i.e., useless) when indels are not evenly distributed, that is: 
when an indel occured in an anchor, or two (or more)
indels occured in a specific island over all $d$ children. 
This is a rare event. Barring this event, we expect that each
island suffers only at most one indel over all children. The island
(of a child) that has exactly one indel is invalidated (second mode of
failure), and is thus deemed useless for reconstruction purposes. As
long as the parent node is not invalidated, each island will have at
least $d-2$ non-invalidated children islands (one additional island is
potentially lost to a child node that may have been invalidated
at an earlier stage).

%\cnote{Even if a child island has only one indel it still adds noise to the reconstruction. ???}

Even when the algorithm identifies that a child island has an indel somewhere, the island is not ignored. The algorithm still needs to compute the length of the island in order to know the start of the next island in this child. For this purpose, we use the anchor of the next island and match it to the
corresponding anchors of the other (non-invalidated) child islands. In
fact the same procedure lets us detect which of the child islands are
invalidated.

More formally, we define $d$ functions
$f_i :\{1, \ldots, \kref_0\} \rightarrow \{1, \ldots, \kref_i\} \cup \{\dagger\}$,
where $f_i$ takes the sites of $v$ to the corresponding sites of the $i$'th child
or to the special symbol $\dagger$ if the site was deleted.  Note that
for each $i$, $f_i$ is monotone, when ignoring sites which are mapped
to $\dagger$.
%\srnote{To emphasize that we compute shifts only for islands, I'm
%  making this a function of $r$ instead of $t$.}
For $t = \lisle r$,
let $s_i(r) = f_{i}(t+1) - (t+1)$ denote the displacement of the site
corresponding to the $(t+1)^{\mathrm{st}}$ site of the parent, in the $i^{\mathrm{th}}$
child.  By convention, we take $s_i(0) = 0$.
%\srnote{I think we need
%  some convention for $0$.}
If there is no indel between $t = \lisle
r$ and $t' = \lisle r'$ then $s_i(r) = s_i(r')$. Note that, in the
specific case of one indel operation in the island, we have that
$|s_i(r) - s_i(r')|=1$.

\noindent\textbf{Algorithm.} 
{Our algorithm estimates the values of $s_i(r)$ and uses these
estimates to match the starting positions of the islands in the
children.}
The full algorithm is given in
Figure~\ref{fig:algorithm} in the Appendix. We use the following additional notation.
For $x\in\{0,1\}$ we let $\spin{x} = 2x - 1$.
Then, for two $\{0,1\}^{m}$-sequences $Y = y_1,\ldots,y_m$ and $Z =
z_1,\ldots,z_m$, we define their (empirical) correlation as
\begin{equation*}
\cor(Y,Z) = \frac{1}{m}\sum_{j=1}^m \spin{y_j}\spin{z_j}.
\end{equation*}
Note that $y \mapsto \spin{y}$ maps $1$ to $1$ and $0$ to $-1$.
One can think of $\cor(Y,Z)$ as a form of normalized centered Hamming distance between
$Y$ and $Z$.
In particular, a large value of $\cor(Y,Z)$ implies that $Y$ and $Z$ tend to agree.
We will use the following threshold (which will be justified in Section~\ref{sec:anchor})
\begin{equation*}
\thres = ((1 - \delta)(1-2\prbs)^2 - 4 \beta),
\end{equation*}
where $\delta$ is chosen so that
$$(1 - \delta)(1-2\prbs)^2 - 8\beta > \delta + 8\beta,$$
where again $\beta = O(d^{-1} \log d)$.
%\aanote{can we give an estimate of $\beta$ here? in fact, in the
%  algorithm, it would be good to give estimates of all quantities,
%  like $a$. alternatively, the main theorem says something like ``for
%  this choices of $a,\beta,...$, the algorithm works''.}

%\section{Analysis of the Reconstruction Process}\label{sec:analysis}
\section{Analyzing the Indel Process}\label{sec:analysis}

%\aanote{Same proof overview as in the intro?}

%Recall that $\kref$ is the number of sites in the root.
We define $a \geq
C \log\nref$
and $\alpha \leq \eps/d < 1$,
for constants $C,\eps$ to
be determined later.
%\srnote{Why is $\alpha$ a function of $1/d$?}
We require $a < \kref^{1/3} < \poly(n)$.\footnote{A variant
of the algorithm where the anchors have length $O(\log k)$
 also works when $k \gg n$.}
%\srnote{Why do we require $\kref <
%  n$?}
We assume that the indel probability per site satisfies
$$\prbid = \frac{\alpha}{4 d \kref^{2/3} a} = O\left(
\frac{1}{\kref^{2/3} \log\nref} \right).$$ Throughout, we denote the
tree by $T = (V,E)$.

\subsection{Bound on the Sequence Length}

As the indel probability is defined per site, longer sequences
suffer more indel operations than shorter ones. We begin by bounding
the effect of this process. We show that with high probability the
lengths of all sequences are roughly equal.
\begin{lemma}[Bound on sequence length]\label{min-max-length}
For all $\zeta > 0$ (small), there exists $C' > 0$ (large)
so that for all $u$ in $V$, we have
\begin{equation*}
k_v \in [\mink,\maxk] \equiv [(1-\zeta)\kref,(1+\zeta)\kref],
\end{equation*}
with high probability given
$\kref \geq C' \log^3\nref$.
We denote this event by $\lcal$.
%\srnote{Covering all cases for $\prbid$ with Chernoff is a pain. A constant approximation is enough here.}
\end{lemma}
%\begin{proof}
%See Appendix.
%\end{proof}

%We condition on the event that all vertices have roughly the same
%length (up to $\sqrt{k}$). As the probability of this event is high, we
%will ignore this in the rest of the analysis.\aanote{do we still need
%  this paragraph?}

\subsection{Existence of a Dense Stable Subtree}

In this section, we show that with probability close to $1$ there exists
a dense subtree of $T$ with a ``good indel structure,'' as defined below. 
Our algorithm will try to identify this subtree and perform reconstruction on it, as described in Section~\ref{sec:ancestral}.

\noindent\textbf{Indel structure of a node.}
Recall that $\lisle = \kref^{1/3}$.
\begin{definition}[Indel structure]
For a node (parent) $v$, we say that $v$ is radioactive if one of the following events happen:
%\srnote{I removed the formal
%  statements? I find them more confusing than anything (eps. the third one). See file.}
\begin{enumerate}
  \item Event $\bcal_1$: Node $v$ has a child $u$ such that when evolving from $v$ to $u$
    an indel operation occurred in at least one of the sites which are
    located in an anchor.
    %Formally, there is some integer $r$ such that
    %$f_u(r\lisle)=\dagger$ or $f_u([r\lisle:r\lisle+a-1])\neq f_u(r\lisle)+\{0,\ldots a-1\}$.
%% $f_{U}(r) - f_{U}(t) \neq r-t$, where $c
%%     m^{1/3} \leq r,t \leq c m^{1/3} + a$, and $c$ is any integer.

  \item Event $\bcal_2$: There is an island $I$ and two children $u,u'$, such
    that an indel occurred in $I$ in the transition from $v$ to $u$
    and in the transition from $v$ to $u'$.
    %Formally, there
    %exists an integer $r$ such that $f_u([r\lisle:(r+1)\lisle-1])\neq
    %f_u(r\lisle)+\{0,\ldots \lisle-1\}$ and $f_{u'}([r\lisle:(r+1)\lisle-1])\neq
    %f_{u'}(r\lisle)+\{0,\ldots \lisle-1\}$.

%%       \[f_{U}(r) - f_{U}(t) \neq r-t\]
%%       \[f_{\tilde U}(\tilde r) - f_{\tilde U}(\tilde t) \neq \tilde r - \tilde t\]
%%       and $c m^{1/3} \leq r,t,\tilde r, \tilde t < (c+1)m^{1/3}$

   \item Event $\bcal_3$: There is an island $I$ and a child $u$, such that two indel
     operations (or more) happened in $I$ in the transition from $v$ to
     $u$.
     %Formally there exist three indices $r,s,t$ such that
     %$f_{U}(r) - r \neq f_{U}(s) - s$, $f_{U}(t) - t \neq f_{U}(s) -
     %s$, $f_{U}(r) - r \neq f_{U}(t) - t$ with $cm^{1/3} \leq r,s,t, <
     %(c + 1)m^{1/3}$ where $c$ is some integer.

\end{enumerate}
Otherwise the node $v$ is stable.
By definition, the leaves of $T$ are stable.
A subtree of $T$ is stable if all of its nodes are stable.
\end{definition}

\begin{lemma}[Bound on radioactivity]
\label{lem:badIndelStructure}
For all $\alpha > 0$,
there exists a choice of $\zeta > 0$ small enough in Lemma~\ref{min-max-length} such that conditioning on the event $\lcal$ occuring:
any vertex $v$ is radioactive with probability at most
$\alpha$.
%\srnote{I think we can leave the comment about independence of substitutions
%to where we use it in the proof. The previous statement made no sense as a probabilistic claim.}
%, where the probability is taken only on the indel structure
%of the node $v$.
%This probability is independent from the actual
%value of $X_v$ and the substitutions.
\end{lemma}
%\begin{proof}
%See Appendix.
%\end{proof}
As a corollary we obtain the following.
%We first need a definition.
%\begin{definition}[Stable Nodes]
%A node $v$ is \emph{stable} if:
%\begin{itemize}
%\item it has a good indel structure; and
%\item it has at most 2 children who are not stable.
%\end{itemize}
%\end{definition}
\begin{lemma}[Existence of a dense stable subtree]\label{lem:stable}
For all $\chi > 0$, there is of $\zeta > 0$ small enough in Lemma~\ref{min-max-length} such that, conditioning on the event $\lcal$ occuring, with probability at least $1-\chi$,
the root of $T$ is the father of
a $(d-1)$-ary stable subtree of $T$. We denote this event by $\scal$.
\end{lemma}
%\begin{proof}
%See Appendix.
%\end{proof}

%For $\alpha$ going to $0$, see Appendix.

%\subsection{Recursive Majority Against an Adversary}\label{sec:ancestral}

\section{A Stylized Reconstruction Process}\label{sec:ancestral}

In this subsection, we lay out the basic lemmas that we need to analyze our
ancestral reconstruction method. We do this by way of describing a hypothetical sequence reconstruction process performed on the stable tree defined by the indel process (see Lemma~\ref{lem:stable}). We analyze this reconstruction process (assuming that the radioactive nodes and the islands with indels are controlled by an adversary) and show in Lemma~\ref{lem:bernoulli} that the process gives strong reconstruction guarantees. Then we argue in Section~\ref{sec:true reconstruction} that our algorithm performs at least as well as the reconstruction process against the adversary described in this section. Throughout this section we suppose that a stable tree exists and is given to us, together with the ``orbit'' of every site of the sequence at the root of the tree (see function $F$ below). However, we are given no information about the substitution process.

Let $v \in V$ and assume $v$ is the root of
a $(d-1)$-ary stable subtree $T^* = (V^*,E^*)$ of $T$.
(We make the stable subtree below $v$ into a $(d-1)$-ary tree
by potentially removing arbitrary nodes from it, at random.) Let $u \in V^*$. 
For each island $I$ in $u$, at most
one child $u'$ of $u$ in $T^*$ contains an indel in which case it contains exactly one indel. 
We say that such an $I$ is a corrupted island of $u'$. The basic intuition
behind our analysis is that, provided the alignment on $T^*$ is performed
correctly (which we defer to Section~\ref{sec:correctness}), the ancestral reconstruction
step of our algorithm is a recursive majority procedure against an adversary
which controls the corrupted islands and the radioactive nodes (as well as all their
descendants). Below we analyze this adversarial process.

\noindent\textbf{Recursive majority.} We begin with a formal definition
of recursive majority. Let $\Maj :
\{0,1,\sharp\}^d \to \{0,1\}$ be the function that returns the majority value
over non-$\sharp$ values, and flips an unbiased coin in case of a tie (including
the all-$\sharp$ vector).
Let $\nref_0 = d^{\Href_0}$ be the number of leaves in $T$ below $v$.
Consider the following
recursive function of $z = (z_1,z_2,\ldots,z_{\nref_0}) \in \{0,1,\sharp\}$:
$\Maj^0(z_1) = z_1$,
and
\begin{equation*}
\Maj^j(z_1,\ldots,z_{d^{j}})
= \Maj(
\Maj^{j-1}(z_1,\ldots,z_{d^{(j-1)}}),
\ldots,
\Maj^{j-1}(z_{d^{j} - d^{(j-1)} + 1},\ldots,z_{d^{j}})),
\end{equation*}
for all $j = 1,\ldots,{\Href_0}$.
Then, $\Maj^{\Href_0}(z)$ is the $d$-wise recursive majority of $z$.

Let $\X_0 = \x^0_1,\ldots,\x^0_{\kref_0}$ be the sequence at $v$. For $u\in V^*$
and $t = 1,\ldots,\kref_0$, we denote by $F_u(t)$ the position of site $\x_t^0$
in $u$ or $\dagger$ if the site has been deleted on the path to $u$.
We say that $\ccal_{u,t}$
holds if $F_u(t)$ is in a corrupted island of $u$.
Let $\path(u,v)$ be the set of nodes on the path between $u$ and $v$.
\begin{definition}[Gateway node]
A node $u$ is a \emph{gateway} for site $t$ if:
\begin{enumerate}
\item $F_u(t) \neq \dagger$; and
\item For all $u' \in \path(u,v) - \{v\}$, $\ccal_{u^{\prime},t}$ does not hold.
\end{enumerate}
\end{definition}
We let $T^{**}_t = (V^{**}_t,E^{**}_t)$ be the subtree of $T^*$ containing
all gateway nodes for $t$. By construction, $T^{**}_t$ is at least $(d-2)$-ary
and for convenience we remove arbitrary nodes, at random, to make it exactly $(d-2)$-ary.
Notice that, for $t,t' \in [1:\kref_0]$, the subtrees $T^{**}_{t}$ and $T^{**}_{t'}$
are random and correlated. However, they are independent of the substitution process.

We will show in Section~\ref{sec:correctness} that the reconstructed sequence produced by our method
at $v$ ``dominates'' (see below) the following reconstruction process. Let
$L_v = u_1,\ldots,u_{\nref_0}$ be the leaves below $v$ ordered according to
a planar realization of the subtree below $v$. Denote
by $\X_{i} = x^i_1,\ldots,x^i_{\kref_i}$ the sequence at $u_i$.
For $t= 1,\ldots,\kref_0$,
let $L^{**}_t$ be the leaves of $T^{**}_t$. We define
the following auxiliary sequences:
for $u_i\in L_v$, we let $\Xi_i = \xi^i_1,\ldots,\xi^i_{\kref_i}$ where
for $t=1,\ldots,\kref_0$
\begin{displaymath}
\xi^i_t =
\left\{
\begin{array}{ll}
\x^i_{F_{u_i}(t)} & \mbox{if $u_i \in L^{**}_t$}\\
1 - \x^0_t & \mbox{o.w.}
\end{array}
\right.
\end{displaymath}
In words, $\xi^i_t$ is the descendant of $\x^0_t$ if $u_i$ is a gateway to $t$
and is the opposite of the value $\x^0_t$ otherwise. Because of the monotonicity
of recursive majority, the latter choice is in some sense the ``worst adversary''
(ignoring correlations between sites---we will come back to this point later).
We then define a reconstructed sequence at $v$ as $\RXi_0 = \rxi^0_1,\ldots,\rxi^0_{\kref_0}$
where for $t=1,\ldots,\kref_0$
\begin{equation*}
\rxi^0_t = \Maj^{\Href_0}(\xi^1_t,\ldots,\xi^{\nref_0}_t).
\end{equation*}
We now analyze the accuracy of this (hypothetical) estimator---which
we refer to as the \emph{adversarial reconstruction of $X_0$}.
We show in Section~\ref{sec:correctness} that our actual
estimator is at least as good as $\RXi_0$ w.h.p.

%\noindent\textbf{Preliminary Claims.}

%\begin{claim}\label{claim:alpha}
%For all $\alpha > 0$, there exists a $\prbid = \prbid(\alpha) > 0$
%small enough such that with probability at least $1-\alpha$
%the indel-free subtree of $\tree{d}{\ell}$ contains a copy of $\tree{d-1}{\ell}$.
%In that case, we say that \emph{the indel-free tree is dense.}
%\end{claim}
%\begin{proof}
%Follows from lemmas in \cite{Mossel:01}.
%\end{proof}

\subsection{Recursive Majority Against an Adversary}

%We analyze the performance of the recursive majority function applied to our leaf sequences, when---being pessimistic---an adversary controls the radioactive nodes and the corrupted islands, as described above. 
To analyze the performance of the adversarial reconstuction $\RXi_0$, we consider the following stylized process.
\begin{definition}[Recursive Majority Against an Adversary] \label{def:recursive majority adversary}
We consider the following process:
\begin{enumerate}
\item Run the evolutionary process on $\tree{d-2}{\Href_0}$
at one position only starting with root state $0$
without indels, that is, taking $\prbid = 0$.

\item Then complete $\tree{d-2}{\Href_0}$ into
$\tree{d}{\Href_0}$ and associate to each additional node the state $1$.

\item Let
$\arm{d}{\Href_0}$ be the random variable in $\{0,1\}$ obtained
by running recursive majority on the leaf states obtained above.
\end{enumerate}
\end{definition}
We call this process the \emph{recursive majority against an adversary on $\tree{d}{\Href_0}$}. We show the following.
\begin{lemma}[Accuracy of recursive majority]\label{lem:beta}
For all $\beta > 0$, there exists a constant $C'' > 0$ such
that taking
\begin{equation*}
\ths^2 > \frac{C'' \log d}{d},
\end{equation*}
and $d$ large enough, then the
probability that the recursive
majority against an adversary on $\tree{d}{\Href_0}$ correctly reconstructs root state $0$
is at least $1-\beta$ uniformly in $\Href_0$. In comparison,
note that the Kesten-Stigum bound for binary symmetric channels
on $d$-ary trees is $\theta^2 > d^{-1}$~\cite{KestenStigum:66,Higuchi:77}.
\end{lemma}
%\begin{proof}
%See Appendix.
%\end{proof}
As a corollary of Lemma~\ref{lem:beta}, we have the following.
\begin{definition}[Bernoulli sequence]
For $q > 0$ and $m \in \nintgr$,
the $(q,m)$-Ber\-noul\-li sequence
is the product distribution on $\{0,1\}^{m}$ such
that each position is $1$ independently with probability $1-q$.
We denote by $\bern{q}{m}$ the corresponding random variable.
\end{definition}
\begin{lemma}[Subsequence reconstruction]\label{lem:bernoulli}
Assume $v$ is the root of a $(d-1)$-ary stable subtree.
For all $\beta > 0$, choosing $C'' > 0$ as in Lemma~\ref{lem:beta} is such
that the following holds for $d$ large enough.
For $t,m \in \{1,\ldots,\kref_0\}$, let $\Lambda = (\lambda_1,\ldots,\lambda_{m})$ be the
\emph{agreement vector} between the $\RXi_0[t+1:t+m]$ and $\X_0[t+1:t+m]$, that is,
$\lambda_i = 1$ if recursive majority correctly reconstructs position $i$.
Then there is $0 \leq \beta' \leq \beta$ such that $\Lambda\sim \bern{\beta'}{m}$.
(Here, $\beta'$ may depend on $\Href_0$ but $\beta$ does not.)
\end{lemma}
%\begin{proof}
%See Appendix.
%\end{proof}

\subsection{Stochastic Domination and Correlation}\label{sec:estimation}

In our discussion so far we have assumed that a stable tree exists and is given to us, together with the the function $F$. This allowed us to define the stylized recursive majority process against an adversary (Definition~\ref{def:recursive majority adversary}), for which we established strong reconstruction guarantees (Lemmas~\ref{lem:beta} and~\ref{lem:bernoulli}). In reality, we have no access to the stable tree. We are going to construct it recursively from the leaves towards the root. At the same time we will align sequences, discover corrupted islands, and reconstruct sequences of internal nodes. The stylized recursive majority process will be used to provide a ``lower bound'' on the actual reconstruction process. The notion of ``lower bound'' that is of interest to us is captured by {\em stochastic domination}, which we proceed to define formally.

%\smallskip \noindent\textbf{Stochastic domination.}
%We will show in Section~\ref{sec:correctness} that our
%ancestral reconstruction ``dominates'' the above
%adversarial process in the following sense.
\begin{definition}[Stochastic domination]
Let $\X$,$\Y$ be two random variables in $\{0,1\}^{m}$.  We say that
$\Y$ stochastically dominates $\X$, denoted $\X\preceq \Y$, if there
is a joint random variable $(\tilde X,\tilde Y)$ such that the
marginals satisfy $\X \sim \widetilde{\X}$ and $\Y \sim
\widetilde{\Y}$ and moreover $\prob[\widetilde{\X} \leq
  \widetilde{\Y}] = 1$.
\end{definition}
Note that in the definition above $X$ and $Y$ may (typically) live in different probability spaces.
Then, the joint variable $(\widetilde{\X},\widetilde{\Y})$ is a coupled version of $X$ and $Y$. In our case,
$X$ is the adversarial recursive process whereas $Y$ is the actual reconstruction performed
by the algorithm.
%These are defined separately without any relation to each other,
%but we will show that we can run them ``with the same coin flips'' in order to obtain
%the monotonicity condition above.
We now show how to use this property for correlation estimation.

\smallskip \noindent\textbf{Correlation.}
The analysis of the previous section guarantees that the sequences output by the adversarial reconstruction process are well correlated with the true sequences. But if we are only going to use the adversarial process as a lower bound for the true reconstruction process, it is important to establish that stochastic domination preserves correlation. In preparing the ground for such a claim let us establish an important property of the adversarial process.
Let $T_u$ and $T_v$ be the two disjoint copies of $T^{(d)}_{h}$ rooted
at the nodes $u$ and $v$ respectively, and let $\X=\x_1,\x_2,\ldots,\x_m \in
\{0,1\}^m$ and $\Y=\y_1,\y_2,\ldots,\y_m \in \{0,1\}^m$ be sequences at the
nodes $u$ and $v$. Assume that $u$ and $v$ are the roots of $(d-1)$-ary stable subtrees.
%Suppose that the sequences $X$ and $Y$ evolve using
%the model of sequence evolution (Definition \ref{def: model of evolution})
Let $\RX'= \rx'_1,\rx'_2,\ldots,\rx'_m \in \{0,1\}^m$ and
$\RY'= \ry'_1,\ry'_2,\ldots,\ry'_m \in \{0,1\}^m$ be the
reconstructions of $\X$ and $\Y$ obtained by the adversarial reconstruction process. Let $\Lambda = \lambda_1,\ldots,\lambda_m$ and $\Theta = \theta_1,\ldots,\theta_m$ be the resulting agreement vectors. We show the following:
%\srnote{I'm splitting this lemma into two to allow conditioning on
%what happens to adversarial reconstruction---a high-probability event---from which
%any dominating recontruction (which we haven't analyzed at this point) is bounded deterministically. In the previous
%version, those two things were mixed.}
\begin{lemma}[Concentration of bias]\label{lem:biasprelim}
Let $\beta',\beta$ be as in Lemma~\ref{lem:bernoulli}. Then,
with probability at least $1-e^{-\Omega(m \beta^2 )}$
the following are satisfied
$$\left|\frac{1}{m}\sum_{i=1}^m \spin{\lambda_i} \spin{\theta_i} - (1-2\beta')^2 \right| \le \frac{1}{2} \beta;$$
$$\left|\frac{1}{m}\sum_{i=1}^m \ind_{\spin{\lambda_i} = -1} - \beta' \right| \le \frac{1}{2} \beta;
\qquad\left|\frac{1}{m}\sum_{i=1}^m \ind_{\spin{\theta_i} = -1} - \beta' \right| \le \frac{1}{2} \beta.$$
%We denote this event by $\bcal_{u,v}$.
\end{lemma}
%\begin{proof}
%See Appendix.
%\end{proof}
We use the previous lemma to argue that stochastic domination does not affect our correlation computations.
\begin{lemma}[Correlation bound]\label{lem:bias}
Let $\RX,\RY \in \{0,1\}^m$ be random strings defined on the same probability space as $\RX'$ and $\RY'$.
Denote by $Z$ (resp.~$W$) the agreement vectors of $\RX$ (resp.~$\RY$) with $\X$ (resp.~$\Y$).
Assume that $\Lambda\leq Z$ and $\Theta \leq W$ with
probability $1$, where $\Lambda$ and $\Theta$ are the agreement vectors of $\RX'$ and $\RY'$ with $\X$ and $\Y$ as explained above.
%(Note that this bound depends on the behavior
%of $\RX'$ and $\RY'$ rather than $\X$ and $\Y$.)
%(We say that $\RX$ and $\RY$ dominate $\RX'$ and $\RY'$ respectively.)
Then,{
\begin{equation*}
|\cor(X,Y) - \cor(\RX,\RY)|
\leq 1 - \frac{1}{m}\sum_{i=1}^m (\spin{\lambda_i} \spin{\theta_i} -
\ind_{\spin{\lambda_i}=-1} - \ind_{\spin{\theta_i}=-1}),
\end{equation*}}
with probability $1$. Furthermore, conditioned on the conclusions
of Lemma~\ref{lem:biasprelim}, we have, with probability $1$:
\begin{equation*}
|\cor(X,Y) - \cor(\RX,\RY)|
\leq 8\beta.
\end{equation*}
\end{lemma}

\section{Analyzing the True Reconstruction Process} \label{sec:true reconstruction}

We provide the proof of Theorem~\ref{thm:main}. In Section~\ref{sec:anchor}, we show that, if a stable subtree exists, the adversarial reconstructions of aligned anchors exhibit strong correlation signal, while misaligned anchors exhibit weak signal. This holds true for sequences that stochastically dominate the adversarial reconstructions. We use this property to complete the analysis of our reconstruction method in Section~\ref{sec:correctness}.

\subsection{Anchor Alignment}\label{sec:anchor}
%To prove the above lemma we need the following notations and claims.
Consider a parent $v$ that is stable.
Let $i,j$ be two children with sequences $\X_i = \x^i_1,\ldots,\x^i_{\kref_i}$ and
$\X_j = \x^j_1,\ldots,\x^j_{\kref_j}$.
Let $t = \lisle r$ and consider the following subsequences
(of length $a$) at $i$ and $j$
$$\AA^{i}_{r} = \x^{i}[t + s_i(r)+1:  t + s_i(r) + a], \ \ \textrm{and} \ \
\AA^{j}_{r} = \x^{j}[t + s_j(r)+1:  t + s_j(r) + a].$$
These are related (but not identical)
to the definition of anchors in the algorithm of
Section~\ref{sec:algo}. In particular, note that by definition
$\AA^{i}_{r}$ and $\AA^{j}_{r}$ are always aligned, in the
sense that they correspond to the same subsequence of $v$.
Consider also the following subsequences
$$\DD^{j}_{r} = \x^{j}[t + s_j(r)  : t + s_j(r) + a - 1] \  \textrm{ and } \
\II^{j}_{r} = \x^{j}[t + s_j(r) + 2: t + s_j(r) + a + 1].$$
%Let $r,t$ be such that for all $0 \le y \le
%k$ $f_i(y+r) - f_{i}(r) = y$, and $f_j(y+t) - f_{j}(t) = y$. Denote
%the part of the original $i$'th child at place $f_i(r)$ as $\hat
%X_i(r) = \hat x^{i}[f_i(r): f_i(r) + k]$, the reconstructed version
%$X_i(r) = x^{i}[f_i(r):f_i(r) + k]$. Similarly for the $j$'th child
%$\hat X_j(t) = \hat x^{j}[f_j(t): f_j(t) + k]$ and for the
%reconstructed version of the $j$'th child $X_j(t) = x^{j}[f_j(t):
%  f_j(t) + k]$.
These are the one-site shifted subsequences for $j$.
The following lemma bounds the correlation between
these strings. More precisely, we show that
$\AA^{i}_{r}$ is always significantly more correlated
to its aligned brother $\AA^{j}_{r}$ than to the misaligned
ones $\DD^{j}_{r}$ and $\II^{j}_{r}$. This follows from the fact
that the misaligned subsequences are sitewise independent.
%\srnote{Are they actually mutually independent? Is there an easy way to prove it?
%For our purposes it's enough that it's a martingale with respect to a product space.
%I ended up using the method of bounded differences although it seems
%we're losing in the constants with that.}
\begin{lemma}[Anchor correlations]\label{lem:anchor}
For all $\delta > 0$ such that $(1 - \delta)(1-2\prbs)^2 - 8\beta > \delta + 8\beta,$
there is $C > 0$ large enough so that
with $a = C\log \nref$,
the following hold:
\begin{enumerate}
\item \textbf{Aligned anchors.}
%\begin{equation*}
$ \prob\left[\cor(\AA^{i}_{r},\AA^{j}_{r})
> (1 - \delta)(1-2 \prbs)^2 \right]
> 1 - \exp\left(-\Omega(a)\right) = 1 - 1/\poly(\nref).
$
%\end{equation*}
%    Moreover, if $i,j$ are not radioactive
%    \[\Pr(\cor(X^i_{r}, X^j_{r}) < ((1 - \delta)(2p_s - 1)^2 - 4 \beta) k) < 2^{-O \left( \delta \sqrt{k} \right)} + 2^{-O \left( \beta \sqrt{k} \right)} .\]
%    \ahnote{the second part comes from the bound proven by Costis. There is a square there which i don't understand.}
\item \textbf{Misaligned anchors.}
%\begin{equation*}
$ \prob\left[\cor(\AA^{i}_{r},\DD^{j}_{r})
< \delta \right]
> 1-\exp\left(-\Omega(a)\right) = 1-1/\poly(\nref)$,
%\end{equation*}
and similarly for $\II^{j}_{r}$.

%\[\Pr \left( \cor(\hat X^i_{r}, \hat X^j_{t}) > \delta k \right) < 2^{-O \left( \delta \sqrt{k} \right)}.\]
%    And if $i,j$ are not radioactive we also get
%    \[\Pr(\cor(X^i_{r}, X^j_{r}) > (\delta + 4 \beta) k) < 2^{-O \left( \delta \sqrt{k} \right)} + 2^{-O \left( \beta \sqrt{k} \right)} .\]
\end{enumerate}
We denote by $\acal_{i,j,r}$ the above events and their symmetric counterparts
under $i \leftrightarrow j$.
\end{lemma}
%\begin{proof}
%See Appendix.
%\end{proof}
\begin{lemma}[Anchor correlations: Reconstructed version]\label{lem:ranchor}
Let $\RX_i = (\rx^i_\iota)_{\iota=1}^{\kref_i}$ and
$\RX_j = (\rx^j_\iota)_{\iota=1}^{\kref_j}$ dominate the adversarial reconstructions
$\RX'_i$ and $\RX'_j$ of $\X_i$ and $\X_j$, as defined in
Lemma~\ref{lem:bias}.
Let $\RAA^i_r = \rx^i[t + s_i(r)+1:  t + s_i(r) + a]$
and similarly for all other possibilities $\RAA\leftrightarrow \RDD,\RII$ and/or
$i\leftrightarrow j$.
Denote by $\bcal_{i,j,r}$ the event that
the conclusions of Lemma~\ref{lem:biasprelim} hold for $\RX'_i$ and $\RX'_j$
over all pairs of intervals involving
$[t + s_i(r):  t + s_i(r) + a-1]$, $[t + s_i(r)+1:  t + s_i(r) + a]$, and $[t + s_i(r)+2:  t + s_i(r) + a+1]$, with $i \leftrightarrow j$ as necessary.
%$[t:t+a-1]$, $[t+1:t+a]$, and $[t+2:t+a+1]$.
Then, conditioned
on $\bcal_{i,j,r}$ and $\acal_{i,j,r}$ we have
$$\cor(\RAA^{i}_{r},\RAA^{j}_{r})
> (1 - \delta)(1-2 \prbs)^2 - 8\beta,$$
$$\cor(\RAA^{i}_{r},\RDD^{j}_{r})
< \delta + 8\beta,%$$
%and
%$$
\qquad\cor(\RAA^{i}_{r},\RII^{j}_{r})
< \delta + 8\beta,$$
as well as their symmetric counterparts under $i\leftrightarrow j$.
\end{lemma}
%\begin{proof}
%See Appendix.
%\end{proof}

\subsection{Proof of Correctness}\label{sec:correctness}

{
%We conclude the proof of Theorem~\ref{thm:main}, showing that our recursive procedure reconstructs the desired sequence at the root of the tree. 
We show that our recursive procedure reconstructs the desired sequence at the root of the tree whenever a collection of good events occurs. Recall the definitions of the events $\lcal$, $\scal$, $\bcal_{i,j,r}$, $\acal_{i,j,r}$ from Lemmas~\ref{min-max-length},~\ref{lem:stable},%~\ref{lem:bias}, 
~\ref{lem:anchor} and~\ref{lem:ranchor}.~\footnote{Event $\lcal$ guarantees that there is no big variance in the nodes' sequence lengths; event $\scal$ guarantees that a stable $(d-1)$-ary subtree exists; the events $\bcal_{i,j,r}$ guarantee that the adversarial reconstruction process is successful, also in preserving correlations between sequences of nodes; and the events $\acal_{i,j,r}$ guarantee that aligned anchors (across sequences of a node's children) exhibit strong correlation signal, while misaligned anchors give weak correlation signal.} Conditioning on $\lcal$ and
$\scal$, denote by $T^* = (V^*,E^*)$ the stable $(d-1)$-ary subtree
of $T$. Then, for all $v \in V^*$, all pairs of children $i,j$ of $v$
in $T^*$, and all $r = 1,\ldots,\maxk/\ell$,
%\srnote{This contortion is necessary to
%  avoid $\acal_{i,j,r}$ being a ``random'' event.}
we condition on the events $\bcal_{i,j,r}$ and $\acal_{i,j,r}$.  Note that having conditioned
on $\lcal$ there is only a polynomial number of such events, since all sequence lengths are bounded by $\maxk$. (If $r \cdot \ell$ is larger than a node's sequence length we assume that the corresponding events are vacuously satisfied.) Finally recall that, conditioning on $\lcal$, the event $\scal$ occurs with probability $1 - \chi$ and all other events occur with high probability.  We denote the collection of events by $\ecal$.

Conditioning on $\ecal$, the proof of correctness of the algorithm follows from a bottom-up induction. The gist of the argument is the following. Suppose that at a recursive step of the algorithm we have reconstructed sequences for all children of a node $v$, which are strongly correlated with the true sequences (in the sense of dominating the corresponding adversarial reconstructions). Having conditioned on the events $\acal_{i,j,r}$ and $\bcal_{i,j,r}$, it follows then that the correct alignments of anchors exhibit strong correlation signal while the incorrect alignments weak correlation signal. Hence, our correlation tests between anchors discover the corrupted islands and do the anchor alignments correctly (at least for all nodes lying inside the stable tree). Hence the shift functions $\rs_i$'s are correctly inferred, and the reconstruction of $v$'s sequence can be shown to dominate the corresponding adversarial reconstruction. The complete proof details are given in App
 endix~\ref{sec: full proof}.

}

%%We first identify the
%%event in which a node can fail the reconstruction.
%%Then we prove
%%that, conditioned on the bad event not happening, the reconstruction
%%is correct.
%
%
%\begin{lemma}[Main Invariant]
%\label{lem:mainInvariant}
%Let $\gamma=k^{-\Omega(1)}$.  For each node $v$ in the tree, with
%probability at least $1-\gamma$, the reconstruction is good, meaning
%that:
%\begin{itemize}
%\item
%$m_v=\hat m_v$;
%\item
%$X_v\oplus \hat X_v$ is statistically dominated by an independent
%  $(\beta, m)$-Bernoulli sequence.
%\end{itemize}
%
%Moreover, with high probability, for every two nodes that do not have
%ancestral relationship, the events are independent.
%\end{lemma}

%\srnote{Do we need a final proof for the main theorem?}

%\section{Concluding Remarks}
%
%\srnote{Discuss possible extensions.}

%\ahnote{We can take this out if there is no room}
%\section{Concluding Remarks}
%
%We briefly mention a few research directions and extensions. Improving the dependency of the indel probability in $k$ %from $1/k^{2/3}$ to a better $1 / \poly(k)$ can be done by handling a constant number of deletions in a specific %island, and by making the islands bigger (to reduce the number of anchors). Handling small $d$ can be done if the %substitution probability is also small. We do not know how to cope with small $d$ (say 2) when the deletion %probability is close to the threshold required for reconstruction when there are no indels.

\bibliographystyle{alpha}
\bibliography{thesis}

\appendix
\clearpage
\section{Algorithm}

\begin{figure}[htb]
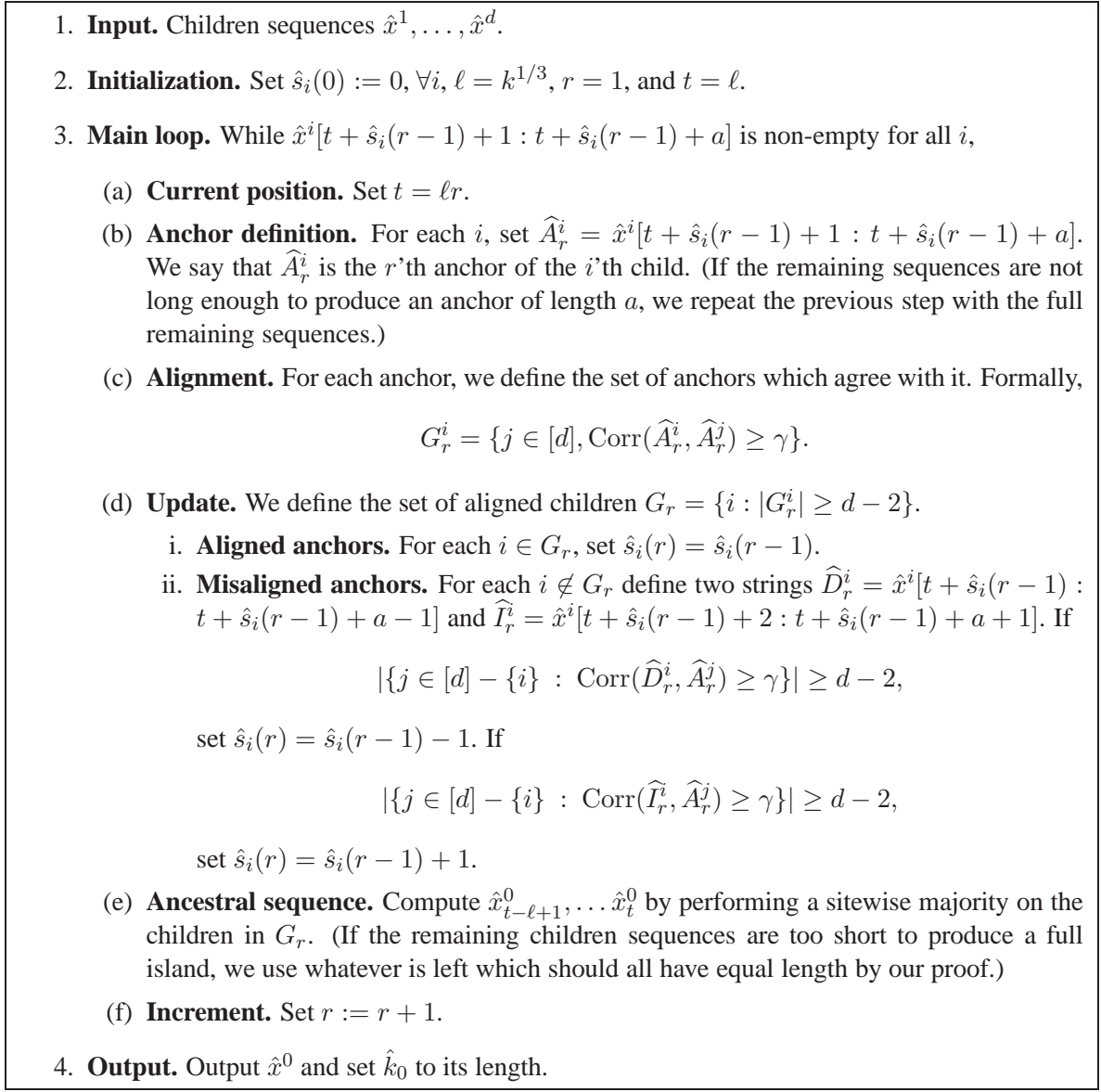

  \centering \fbox{
  \begin{minipage}{6in}
{
\begin{enumerate}
	\item \textbf{Input.} Children sequences $\rx^1,\ldots,\rx^d$.
  \item \textbf{Initialization.} Set $\rs_i(0) := 0$, $\forall i$, $\lisle = k^{1/3}$, $r=1$, and $t=\lisle$.
  \item \textbf{Main loop.} While $\rx^{i}[t + \rs_i(r-1)+1:  t + \rs_i(r-1) + a]$ is non-empty for all $i$,
      \begin{enumerate}
      \item \textbf{Current position.} Set $t = \lisle r$.
        \item \textbf{Anchor definition.} For each $i$, set
        $\RA^{i}_{r} = \rx^{i}[t + \rs_i(r-1)+1:  t + \rs_i(r-1) + a]$.
        We say that $\RA^{i}_{r}$ is the $r$'th anchor of the $i$'th child.
        (If the remaining sequences are not long enough to produce an
        anchor of length $a$, we repeat the previous step with the full remaining sequences.)
        \item \textbf{Alignment.} For each anchor, we define the set of anchors which agree with it.
        Formally, $$G^{i}_{r} = \{j \in [d], \cor(\RA^{i}_{r},\RA^{j}_{r}) \ge \thres\}.$$
        \item \textbf{Update.} We define the set of aligned children $G_{r} = \{i : |G^{i}_{r}| \geq d-2\}$.
				\begin{enumerate}
				\item \textbf{Aligned anchors.} For each $i \in G_r$, set $\rs_i(r) = \rs_i(r-1)$.
        \item \textbf{Misaligned anchors.} For each $i \not \in G_{r}$ define
          two strings $\RD^{i}_{r} = \rx^{i}[t + \rs_i(r-1)  : t +
            \rs_i(r-1) + a - 1]$ and $\RI^{i}_{r} = \rx^{i}[t +
            \rs_i(r-1) + 2: t + \rs_i(r-1) + a + 1]$. If
            \[|\{j \in [d]-\{i\} \ : \ \cor(\RD^{i}_{r},\RA^{j}_{r}) \ge \thres\}| \geq d-2,\]
            set $\rs_i(r) = \rs_{i}(r-1)  - 1$. If
             \[|\{j \in [d]-\{i\} \ : \ \cor(\RI^{i}_{r},\RA^{j}_{r}) \ge \thres \}| \geq d-2,\]
             set $\rs_i(r) = \rs_{i}(r-1)  + 1$.
        \end{enumerate}
        \item \textbf{Ancestral sequence.} Compute $\rx^0_{t - \lisle + 1}, \ldots \rx^0_{t}$ by performing a sitewise majority on the children in $G_r$. (If the remaining children sequences are too short to produce
        a full island, we use whatever is left which should all have equal length by our proof.)
        \item \textbf{Increment.} Set $r := r+1$.
      \end{enumerate}
      \item \textbf{Output.} Output $\rx^0$ and set $\rkref_0$ to its length.
\end{enumerate}
}
  \end{minipage}
  }
  \caption{{This is the basic recursive step of our reconstruction algorithm.
  It takes as input the $d$ inferred sequences of the children $\rx^1, \ldots,
  \rx^d$ and computes a sequence for the parent $\rx^0$. If any of the steps above
  cannot be accomplished, we abort the reconstruction of the parent and declare it
  radioactive.}}
  \label{fig:algorithm}
\end{figure}

\section{Further Lemmas}

For $\alpha$ going to $0$, we have more precisely:
\begin{lemma}[Limit $\alpha \to 0$]\label{lem:limit}
Condition on $\lcal$. Let
\begin{equation*}
\alpha = \frac{1}{h(\nref)},
\end{equation*}
for $h(\nref) = \omega(1)$. Then, for $n$ large enough,
the root is the father of a $(d-1)$-ary stable subtree
with probability at least
\begin{equation*}
1 - \chi = 1 - \frac{1}{\sqrt{h(\nref)}}.
\end{equation*}
\end{lemma}
\begin{app-proof}{Lemma~\ref{lem:limit}}
Plugging $\alpha = 1/h(\nref)$ and $\nu = 1 - 1/\sqrt{h(\nref)}$
into the recursion derived in the proof of Lemma~\ref{lem:stable}, we get
\begin{eqnarray*}
(1 - \alpha)g(\nu)
&=& \left(1 - \frac{1}{h(\nref)}\right)\Bigg(1 - \frac{d}{\sqrt{h(\nref)}} + O\left(\frac{1}{h(\nref)}\right)\\
&& \qquad + d\left(1 - \frac{d-1}{\sqrt{h(\nref)}} + O\left(\frac{1}{h(\nref)}\right)\right)\frac{1}{\sqrt{h(\nref)}}\Bigg)\\
&=& \left(1 - \frac{1}{h(\nref)}\right) \left(1 - O\left(\frac{1}{h(\nref)}\right)\right)\\
&\geq& 1 - 1/\sqrt{h(\nref)},
\end{eqnarray*}
for $n \to +\infty$.
\end{app-proof}

\section{Proofs}

\begin{app-proof}{Lemma~\ref{min-max-length}}
We prove the upper bound by assuming there is no deletion. The lower bound
can be proved similarly. The proof goes by induction. Let $v$ be a node
at graph distance $i$ from the root. We show that there is $C'' > 0$ independent of $i$
such that
\begin{equation*}
k_v \leq \kref + i \sqrt{C'' \kref \log\nref}.
\end{equation*}
Since the depth of $T$ is $O(\log\nref)$,
this implies the main claim as long as
\begin{equation*}
\sqrt{C'' \kref \log\nref} \log\nref \leq \zeta \kref,
\end{equation*}
which follows from our assumption for $C' > 0$ large enough.

The base case of the induction is satisfied trivially.
Assume the induction claim holds for $v$, the parent of $u$.
It suffices to show that the number of new insertions is at most
$\sqrt{C'' \kref \log\nref}$. By our induction hypothesis, the number of insertions
is bounded above by a binomial $Z$ with parameters $\kref + (i-1) \sqrt{C'' \kref \log\nref}
\leq (1+\zeta) \kref$
and $\prbid$ w.h.p. By Hoeffding's inequality,
taking
\begin{equation*}
\eta = \sqrt{\frac{C''' \log\nref}{(1+\zeta) \kref}},
\end{equation*}
we have
\begin{eqnarray*}
\prob[Z > (1+\zeta)\kref\prbid + (1+\zeta)\kref\eta]
&<& \exp(-2((1+\zeta)\kref\eta)^2/[(1+\zeta)\kref])\\
&=& 1/\poly(\nref).
\end{eqnarray*}
By our assumption on $\prbid$, we have
\begin{equation*}
(1+\zeta)\kref\prbid = O\left(\frac{\alpha \kref^{1/3}}{\log\nref}\right),
\end{equation*}
so that choosing $C''$ large enough gives
\begin{equation*}
(1+\zeta)\kref\prbid + (1+\zeta)\kref\eta \leq \sqrt{C'' \kref \log\nref}.
\end{equation*}
This proves the claim.
\end{app-proof}

\begin{app-proof}{Lemma~\ref{lem:badIndelStructure}}
According to Lemma~\ref{min-max-length},
the length of the sequence at $v$
is in $[\mink,\maxk]$ w.h.p.
We denote that event by $\lcal_v$.
% Let also $\bcal_1, \bcal_2, \bcal_3$ be the three events leading to a radioactive $v$.
We bound the probability of events $\bcal_1, \bcal_2, \bcal_3$ separately.

Let $\maxnisle = \maxk/\ell = (1+\zeta)\kref^{2/3}$.
Conditioned on $\lcal_v$, there are at most
$\maxnisle$ anchors, each of length $a$. By a union bound, the
probability that at least one of the sites in the anchors has an
indel operation in any child is upper bounded by
\begin{eqnarray*}
\prob[\bcal_1]
&=& \prob[\bcal_1\,|\,\lcal_v]\prob[\lcal_v] + \prob[\bcal_1\,|\,\lcal_v^c]\prob[\lcal_v^c]\\
&\leq& \maxnisle a d \prbid + 1/\poly(\nref)\\
&=& \frac{\alpha a d \maxnisle}{4\kref^{2/3} a d} + 1/\poly(\nref)\\
&=& \frac{(1+\zeta)\kref^{2/3}}{\kref^{2/3}}\cdot \frac{\alpha}{4} + 1/\poly(\nref)\\
&<& \alpha(1/3 - 1/\poly(\nref)),
\end{eqnarray*}
where we choose $\zeta$ small enough.
The quantity we want to estimate is in fact $\prob[\bcal_1\,|\,\lcal]$
(which is not the same as conditioning on $\lcal_v$ only).
But notice that
\begin{equation*}
\prob[\bcal_1]
= \prob[\bcal_1\,|\,\lcal]\prob[\lcal]
+ \prob[\bcal_1\,|\,\lcal^c]\prob[\lcal^c]
\geq \prob[\bcal_1\,|\,\lcal]\prob[\lcal],
\end{equation*}
which implies
\begin{equation*}
\prob[\bcal_1\,|\,\lcal] \leq \frac{\alpha(1/3 - 1/\poly(\nref))}{1 - 1/\poly(\nref)} < \alpha/3.
\end{equation*}
(This argument shows that it suffices
to condition on $\lcal_v$. We apply the same trick below.)

To bound the probability of the second event, consider an island $I$
and a son $u$. The probability that there is an indel when evolving
from $v$ to $u$ is at most
\begin{equation*}
\prbid \lisle = \frac{\alpha}{4\kref^{2/3} a d} \kref^{1/3} = \frac{\alpha}{4\kref^{1/3} a d}.
\end{equation*}
Thus, the probability that more than one child of $v$ experiences an indel in
$I$ is at most
%\srnote{I couldn't understand your calculation so I changed it a little bit.}
\begin{eqnarray*}
\sum_{i = 2}^{d} {d \choose i} \left(\frac{\alpha}{4\kref^{1/3} a d}\right)^i
&\leq& \sum_{i = 2}^{d} \frac{d^i}{i!} \left(\frac{\alpha}{4\kref^{1/3} ad}\right)^i\\
&\leq& \sum_{i = 2}^{d} \frac{1}{i!} \left(\frac{\alpha}{4\kref^{1/3} a}\right)^i\\
&\leq& e \left(\frac{\alpha}{4\kref^{1/3} a}\right)^2\\
&=& \frac{e \alpha^2}{16\kref^{2/3} a^2},
\end{eqnarray*}
where we used that the expression in parenthesis on the second line is $<1$.
Taking a union bound over all islands, the probability that at least two
children experience an indel in the same island is at most
\begin{eqnarray*}
\prob[\bcal_2\,|\,\lcal]
&\leq&
\maxnisle \cdot \frac{e \alpha^2}{16\kref^{2/3} a^2}\\
&=& \frac{(1+\zeta)e \alpha^2}{16 a^2}\\
&<& \frac{\alpha}{3},
\end{eqnarray*}
where we used that $\alpha < 1$.

For the third event, consider again an island $I$ and a child $u$. The
probability for at least two indel operations in $I$ when evolving
from $v$ to $u$ is at most
%\srnote{I multiplied the number of sites by $2$ to account
%for both insertions and deletions.}
\begin{eqnarray*}
\sum_{i = 2}^{2\lisle} {2\lisle \choose i} \left( \frac{\alpha}{4 a d \kref^{2/3}} \right)^i
&\le& \sum_{i = 2}^{2\lisle} \frac{1}{i!} \left( \frac{2\lisle \alpha}{4 a d \kref^{2/3}} \right)^i\\
&\le& \sum_{i = 2}^{2\lisle} \frac{1}{i!} \left( \frac{\alpha}{2 a d \kref^{1/3}} \right)^i\\
&\le& e \left( \frac{\alpha}{2 a d \kref^{1/3}} \right)^2\\
&\le& \frac{e \alpha^2}{4 a^2 d^2 \kref^{2/3}}.
\end{eqnarray*}
(We use $2 \lisle$ to account for insertions {\em and} deletions.)
Taking a union bound over all islands and children, the
probability that there are two indel operations in the same child in
the same island is bounded by
\begin{eqnarray*}
\prob[\bcal_3\,|\,\lcal]
&\leq& d \maxnisle \frac{e \alpha^2}{4 a^2 d^2 \kref^{2/3}}\\
&\leq& \frac{(1+\zeta) e \alpha^2}{4 a^2 d}\\
&<& \alpha/3.
\end{eqnarray*}
Taking a union bound over the three ways in which a site can become
radioactive proves the lemma.
\end{app-proof}

\begin{app-proof}{Lemma~\ref{lem:stable}}
We follow a proof of~\cite{Mossel:01}.
Let $v$ be a node at distance $r$ from the leaves.
We let $\nu_r$ be the probability that $v$ is the root of
a $(d-1)$-ary stable subtree. Let
\begin{equation*}
g(\nu) = \nu^d + d \nu^{d-1} (1-\nu).
\end{equation*}
Then, from Lemma~\ref{lem:badIndelStructure},
\begin{equation*}
\nu_r \geq (1-\alpha)g(\nu_{r-1}).
\end{equation*}
Note that
\begin{equation*}
g'(\nu) = d(d-2)\nu^{d-2}(1-\nu).
\end{equation*}
In particular, $g$ is monotone, $g(1) = 1$, and $g'(1) = 0$.
Hence, there is $1-\chi < \nu^* < 1$ such that
\begin{equation*}
g(\nu^*) > \nu^*.
\end{equation*}
Then, taking
$$1-\alpha > \nu^*/g(\nu^*),$$
we have
\begin{equation*}
\nu_r \geq (1-\alpha) g(\nu_{r-1}) \geq \frac{\nu^*}{g(\nu^*)} g(\nu_{r-1}) \geq \nu^* > 1 - \chi,
\end{equation*}
by the induction hypothesis that $\nu_{r-1} \geq \nu^*$. Note in particular that
$\nu_0 = 1 \geq \nu^*$.
\end{app-proof}

\begin{app-proof}{Lemma~\ref{lem:beta}}
Recall that we assume the root state is $0$ and all
adversarial nodes are $1$. Because of the bias towards
$1$, we cannot apply standard results about recursive majority
for symmetric channels~\cite{Mossel:98,Mossel:04a}. Instead, we
perform a tailored analysis
of this particular channel.

We take asymptotics as $d \to +\infty$ and we show that
the probability of reconstruction can be taken to be
\begin{equation*}
1 - \beta = 1 - \frac{1}{d},
\end{equation*}
for $C''$ large enough.
%We can restrict ourselves
%to the $(d-2)$-ary subtree of non-adversarial nodes
%which we denote by $\tree{d-2}{\Href_0}$ but
%we require a ``strict majority'' for success. More precisely,
Let $v$ be the root of $\tree{d}{\Href_0}$.
We denote by $Z_v$ the number of non-adversarial children of $v$ in
state $0$ and by $Z'_v$ the number of nodes among them
that return $0$ upon applying recursive majority
to their respective subtree. Let $q^0_{\Href_0}$
be the probability of incorrect reconstruction at $v$ (given that the state at $v$ is $0$).
Then
\begin{eqnarray}
1-q^0_{\Href_0}
&\geq& \prob\left[Z'_v \geq \frac{d+1}{2}\right]\nonumber\\
&\geq& \sum_{i=0}^{d-2} \prob\left[Z'_v \geq \frac{d+1}{2}\,|\,Z_v = i\right]\prob[Z_v = i],\label{eq:q0}
\end{eqnarray}
where we simply ignored the contribution of the children who flipped to $1$.

We prove $q^0_{\Href_0} \leq 1/d$ by induction on the height.
Let
$u$ be a non-adversarial node in $\tree{d}{\Href_0}$ at height
$h$ from the leaves to which we associate as above the variables $Z_u, Z'_u$
and the quantity $q^0_{h}$. Note that $q^0_0 = 0$. We assume the induction
hypothesis holds for $h-1$. Note that conditioned on the state at $u$ being $0$ 
$Z_u$ is $\bin(d-2,(1-\prbs))$ where
\begin{equation*}
1 - \prbs = \frac{1+\ths}{2} = \frac{1}{2} + \Theta\left(\sqrt{\frac{\log d}{d}}\right),
\end{equation*}
as $d \to +\infty$. Similarly, given $Z_u = i$, the variable
$Z'_u$ is $\bin(i,1-q^0_{h-1})$. In particular, the quantity
\begin{equation*}
\prob\left[Z'_u \geq \frac{d+1}{2}\,|\,Z_u = i\right],
\end{equation*}
is monotone in $i$.
We use Chernoff's bound on $Z'_u$
to truncate the lower bound (\ref{eq:q0}).
Indeed, let
\begin{equation*}
\mu = (1 - \prbs)(d-2) = \frac{d}{2} + \Upsilon(d),
\end{equation*}
with
\begin{equation*}
\Upsilon(d) = \Theta(\sqrt{d \log d}),
\end{equation*}
and
\begin{equation*}
\mu(1 -\eta) = \frac{d}{2} + \frac{\Upsilon(d)}{2},
\end{equation*}
where in particular
\begin{equation*}
\eta = \Theta\left(\sqrt{\frac{\log d}{d}}\right).
\end{equation*}
Then, we have
\begin{equation*}
\prob[Z_u < \mu(1 -\eta)]
< \exp\left(-\mu \eta^2/2\right) = d^{-\Omega(1)},
\end{equation*}
for $C''$ large enough. Applying to (\ref{eq:q0}) leads to the lower bound
\begin{equation*}
1-q^0_h \geq (1 - d^{-\Omega(1)})\prob\left[\bin\left(\frac{d}{2} + \frac{\Upsilon(d)}{2}, 1-q^0_{h-1}\right) \geq \frac{d+1}{2}\right].
\end{equation*}
By the induction hypothesis, $q^0_{h-1} \leq 1/d$. By applying Chernoff's bound again we get
\begin{equation*}
\prob\left[\bin\left(\frac{d}{2} + \frac{\Upsilon(d)}{2}, 1-q^0_{h-1}\right) \geq \frac{d+1}{2}\right] > 1 - d^{-\Omega(1)},
\end{equation*}
and therefore $q^0_h \leq 1/d$. This proves the claim.
\end{app-proof}

\begin{app-proof}{Lemma~\ref{lem:bernoulli}}
As we pointed out earlier, although the subtrees $(T^{**}_{t'})_{t'=t+1}^{t+m}$
are correlated by the construction of the islands, they are independent
of the substitution process. By forcing (randomly) the subtrees $(T^{**}_{t'})_{t'=t+1}^{t+m}$ to be $(d-2)$-ary and fixing the adversarial nodes to $1$, we restore
the i.i.d.~nature of the sites, from which the result follows.
\end{app-proof}

\begin{app-proof}{Lemma~\ref{lem:biasprelim}}
This follows from Lemma~\ref{lem:bernoulli}, the independence of $\Lambda$ and $\Theta$,
and three applications of Hoeffding's lemma.
\end{app-proof}

\begin{app-proof}{Lemma~\ref{lem:bias}}
Note that
$$\cor(\RX,\RY) = \frac{1}{m} \sum_{i=1}^m \spin{\rx_i} \spin{\ry_i}
= \frac{1}{m}
\sum_{i=1}^m \spin{\x_i} \spin{\y_i} \spin{z_i} \spin{w_i}.$$
Hence,
\begin{equation*}
|\cor(X,Y) - \cor(\RX,\RY)|  \le \frac{1}{m} \sum_{i=1}^m
(1 - \spin{z_i} \spin{w_i}) = 1 -  \frac{1}{m} \sum_{i=1}^m \spin{z_i} \spin{w_i}. %\label{eq: difference bound}
\end{equation*}
Now notice by case analysis that
\begin{equation*}
\spin{z_i} \spin{w_i} \ge \spin{\lambda_i} \spin{\theta_i} -
\ind_{\spin{\lambda_i}=-1} - \ind_{\spin{\theta_i}=-1}. %\label{eq: coupling of bound with Bernoulli process}
\end{equation*}
This proves the first claim.
The second claim follows from the bounds in Lem\-ma~\ref{lem:biasprelim}.
\end{app-proof}

\begin{app-proof}{Lemma~\ref{lem:anchor}}
For the first claim, note that
\begin{equation*}
\expec[\cor(\AA^{i}_{r},\AA^{j}_{r})] = \ths^2 = (1-2 \prbs)^2,
\end{equation*}
where we used that 1) there is no indel in the sites $[t+1: t+a]$
between $v$ and $i,j$; 2) that the sites are perfectly aligned; and 3) that
the substitution process is independent of the indel process.
We also used the well-known fact that the $\ths$'s are
multiplicative along a path under our model
of substitution~\cite{SempleSteel:03}.
The result then follows from Hoeffding's inequality.

For the second claim, because the anchors are now misaligned
the $t'$-th term in $\cor(\AA^{i}_{r},\DD^{j}_{r})$ for $t' \in [t+1: t+a]$
is the variable $\spin{\x^i_{t' + s_i(r)}}\spin{\x^j_{t' + s_j(r) -1}}$
which is uniform in $\{-1,+1\}$.
In particular, we now have
\begin{equation*}
\expec[\cor(\AA^{i}_{r},\DD^{j}_{r})] = 0.
\end{equation*}
The result follows from the method of bounded differences with
respect to the
independent vectors
$$\{(\x^i_{t' + s_i(r)}, \x^j_{t' + s_j(r)})\}_{t'=t}^{t+a}.$$
%If two sites came from a common ancestor, the expected correlation between them is at least $(2p - 1)^2$. The proof follows by applying a Chernoff bound, and Lemma \ref{lemma:weak adversary}.
\end{app-proof}

\begin{app-proof}{Lemma~\ref{lem:ranchor}}
This follows from Lemmas~\ref{lem:bias} and~\ref{lem:anchor} and
the triangle inequality.
\end{app-proof}

\section{Completing the Proof of the Main Theorem}\label{sec: full proof}

Having conditioned on the event $\ecal$, we justify the correctness of our reconstruction method via the following induction. The top level of the induction establishes Theorem~\ref{thm:main}.

\noindent\textbf{Induction hypothesis.}
Consider a parent $v$ in $T^*$; in particular, $v$
is stable. We assume that the following conditions, denoted by $(\star)$, are satisfied:
For all children $i \in [d]$ of $v$ belonging to $T^*$
\begin{enumerate}
\item \textbf{Alignment.}
For all children $i'$ of $i$ with $i' \in T^*$ and all $r=1,\ldots,\maxk/\ell-1$,
\begin{equation}\label{eq:1}
\rs_{i'}(r) = s_{i'}(r).
\end{equation}
(This condition is trivially satisfied for values of $r \ell$ that are larger than the sequence length of $i'$.)
\item \textbf{Reconstruction.}
Moreover, we have $\rkref_i = \kref_i$
%\srnote{This may have
%to be changed if the last island is dropped.}
and for all $t = 1,\ldots,\kref_i$, the following holds:
\begin{quote}
Let $L_i$ be the leaves below $i$ with $\nref_i = |L_i|$.
Let $\Href$ be the level of $v$. Let $L^{**}_t$ be the gateway
leaves for site $t$. For $u\in L^{**}_t$ let $F_u(t)$ be the position
of site $t$ in $u$. Note that $\rx^i_t$ can be written as $\rx^i_t = \Maj^{\Href-1}(z_1,\ldots,z_{n_i})$,
%\end{equation*}
where $z_j$ is either $\sharp$ or $\x^{j}_{\flat_j}$ for an appropriate
function $\flat_j$. Our hypothesis is that
\begin{equation}\label{eq:2}
\forall u\in L^{**}_t,\ \flat_u = F_u(t).
\end{equation}
\end{quote}
In particular,
the ancestral reconstruction $\RX_i$ dominates
the adversarial reconstruction $\RX'_i$.
\end{enumerate}
The base case where $v$ is a leaf is trivially satisfied.

\bigskip \noindent\textbf{Alignment.}
We begin with the correctness of the alignment.
\begin{lemma}[Induction: Alignment]\label{lem:indali}
Assuming $\ecal$ and $(\star)$,
the algorithm infers $s_i$ correctly for all children $i\in[d]$
which are also in $T^*$, that is, (\ref{eq:1}) holds for $v$.
\end{lemma}

\begin{app-proof}{Lemma~\ref{lem:indali}}
%\srnote{Probably need a bit more detail on this proof.
%Explaining more precisely how we use Lemma~\ref{lem:ranchor}
%may be a good idea.}
Let $\Pi$ denote the
set of children of $v$ in $T^*$.
The proof follows by induction on $r$. The base case $r = 0$ is trivial.
Assume correctness for $r - 1$.

If there is no indel in any of the
children $i \in \Pi$ between the sites
$(r - 1)\lisle$ and $r \lisle$ of $v$,
then under $\ecal$, $(\star)$ and Lemma~\ref{lem:ranchor}
we have $\Pi \subseteq G_r$.
In that case, for all $i\in \Pi$
we have $\rs_{i}(r) = \rs_{i}(r -1) = s_i(r-1) = s_i(r)$,
where the second equality is from $(\star)$.

If there is an indel operation in island $r$,
then since $v$ is stable
only one indel operation occurred in one child.
% (and
%there were no indels in the other children).
Denote the child with an indel by $j$.
Assume the indel is a deletion. (The case of the insertion
is handled similarly.)
If $j$ is not in $T^*$
we are back to the previous case. So assume $j$ is in $T^*$.
Again, from $\ecal$, $(\star)$ and Lemma~\ref{lem:ranchor}
the other children in $T^*$ are added to the set $G_r$, and
the shift value will be computed correctly for them.
Moreover by $(\star)$,
for every $i \in \Pi-\{j\}$,
%\srnote{I couldn't understand the previous argument so I changed it a bit.}
%$f_{i}(r \cdot m^{1/3} -
%m^{1/3}) + s_{i}(r \cdot m^{1/3} - m^{1/3}) = f_{j}(r \cdot m^{1/3} -
%m^{1/3}) + s_{j}(r \cdot m^{1/3} - m^{1/3})$.
%If there was a deletion,
%then $f_{i}(r \cdot m^{1/3}) + s_{i}(r \cdot m^{1/3}) = f_{j}(r \cdot
%m^{1/3}) + s_{j}(r \cdot m^{1/3} - m^{1/3}) -1$,
\begin{eqnarray*}
f_i(r\lisle+1)
&=& r\lisle+1 + s_i(r)\\
&=& r\lisle+1 + \rs_i(r)\\
&=& r\lisle+1 + \rs_i(r-1),
\end{eqnarray*}
which is the starting point of $\RA^{i}_{r}$.
Also,
\begin{eqnarray*}
f_j(r\lisle+1)
&=& r\lisle+1 + s_j(r)\\
&=& r\lisle+1 + s_j(r-1) -1\\
&=& r\lisle+1 + \rs_j(r-1) -1\\
&=& r\lisle + \rs_j(r-1),
\end{eqnarray*}
which is the starting point of $\RD^{j}_{r}$.
Thus according to
Lemma~\ref{lem:ranchor} $\RD^{j}_{r}$ matches $\RA^{i}_{r}$ for all $i \in \Pi\cap G_r$.
As there are $d - 2$ children in $\Pi\cap G_r$, we get that the algorithm sets
\begin{equation*}
\rs_j(r)
=  \rs_j(r-1) - 1 = s_j(r-1) - 1 = s_j(r),
\end{equation*}
as required.
Note also that in this case, according to
Lemma~\ref{lem:ranchor} again, $\RA^{j}_r$ does not have high correlation with $\RA^{i}_{r}$ for any $i \in \Pi\cap G_r$, and thus we will consider $\RI^{j}_{r}$ and $\RD^{j}_{r}$. Similarly, $\RI^{j}_r$ does not have high correlation with $\RA^{i}_{r}$ for any $i \in \Pi\cap G_r$, and thus we will not try to set $\rs_{j}(r)$ twice.
%A similar argument applies if there is an insertion in the $j$'th child.
\end{app-proof}

%We are finally ready to prove Lemma~\ref{lem:anchorAlignments}.
%% , namely
%% the correctness of one level of the reconstruction, assuming that
%% the father is not radioactive, and that at most one of the children in
%% the $d-star$ is radioactive. We make no assumptions on the
%% reconstructed values for the radioactive child.
%
%% The lemma is proved in two stages: first we show that the alignment
%% is correct (Lemma~\ref{alignment}), and then that given the alignment,
%% applying majority gives the required guarantee on the sites of the
%% father (Lemma~\ref{reconstruction}).
%
%\begin{lemma}\label{alignment}
%  With probability $1 - 1/k^{-\poly}$, for all $i \in L$, for all
%  $0 \le r \le m^{2/3} - 1$ we have $f_{i}(m^{1/3} \cdot r) = m^{1/3}\cdot r + s_{i}(r \cdot m^{1/3})$
%\end{lemma}
%\begin{proof}
%See Appendix.
%\end{proof}

\noindent\textbf{Ancestral reconstruction.}
We use Lemma~\ref{lem:indali} to prove that the ancestral
reconstruction dominates the adversarial reconstruction.
In the algorithm, we perform a sitewise majority vote over
the children of $v$ in $G_r$ (these are the aligned children---see the description of the algorithm in Figure~\ref{fig:algorithm}). For notational convenience,
we assume that in fact we perform a majority vote over \emph{all} children
but we replace the states of the children outside $G_r$ with $\sharp$.
\begin{lemma}[Induction: Reconstruction]\label{lem:indrec}
Assuming $\ecal$, $(\star)$ and the conclusion of Lemma~\ref{lem:indali},
(\ref{eq:2}) holds for $v$. In particular,
the ancestral reconstruction $\RX_v$ dominates
the adversarial reconstruction $\RX'_v$.
\end{lemma}

\begin{app-proof}{Lemma~\ref{lem:indrec}}
%\srnote{May want to add more details to this proof also.}
The second claim follows from the first one together with
the construction of the adversarial process and the monotonicity
of majority.

As for the first claim, by Lemma~\ref{lem:indali}
for each site of $v$ there are $d-2$ uncorrupted children islands containing
this site such that the children are also in $T^*$.
In particular,
the $d-2$ corresponding sites in the children are
correctly aligned. Moreover, by the induction hypothesis,
each corresponding site in the children satisfy
(\ref{eq:2}). By taking a majority vote over these sites
we get (\ref{eq:2}) for $v$ as well.

A small technical detail is handling the case where the last island has less than $a$ sites, and thus does not contain an anchor. However, in this case, if the father is stable then there are no indel operations at all in the last island, and therefore aligning it according to the previous one gives the right result.
\end{app-proof}
%\begin{proof}
%See Appendix.
%\end{proof}

\end{document}